\documentclass{amsart}
\usepackage{amssymb}
\usepackage[utf8]{inputenc}
\usepackage{color,fullpage}
\usepackage[hidelinks]{hyperref}
\usepackage{scrextend}

\usepackage{tikz}
\usetikzlibrary{automata, positioning}
\usepackage{systeme}
\usepackage{bbm}
\usepackage{float}
\usepackage{enumitem}
\usepackage{listings}
\usepackage{accsupp}

\author{Vincent Ciarcia, Erik Insko}
\title{Limit Cases And Strategy In Chutes \& Ladders}
\newcommand{\ds}{\displaystyle}

\newenvironment{changemargin}[2]{%
\begin{list}{}{%
\setlength{\topsep}{0pt}%
\setlength{\leftmargin}{#1}%
\setlength{\rightmargin}{#2}%
\setlength{\listparindent}{\parindent}%
\setlength{\itemindent}{\parindent}%
\setlength{\parsep}{\parskip}%
}%
\item[]}{\end{list}}

\definecolor{60color}{RGB}{133,180,95}
\definecolor{94color}{RGB}{0,120,230}
\definecolor{99color}{RGB}{100,120,255}

\usepackage{verbatim}

\begin{document}
\begin{center}
\maketitle
\end{center}
\begin{abstract}
We analyze what happens to the average duration of a game of Chutes \& Ladders as the probability of rolling $\delta \in \lbrace 1,2,3,4,5,6\rbrace$ approaches 100\%. We utilize Markov models, and Monte Carlo simulations in Python. We also introduce strategy to the board game by allowing the player to choose whether or not they flip a coin after each die roll where if they get heads they will advance one square and if they get tails they will go back one square. The strategy the player employs to decide when to flip the coin has a significant impact on the average duration of the game. We analyze six different non-trivial strategies.
\end{abstract}

\section{Introduction}\label{Intro}
Chutes \& Ladders is a board game for children ages 3+. Different variations of the game have been found in many different places all over the world. The game originated in ancient India under the name Moksha Patam with a purpose of teaching children basic Hindu philosophy and how life is complicated by virtues and vices. Later, the game found its way to England under the name Snakes \& Ladders in 1892. Eventually, Milton Bradley released the game in its modern form in the US under the name Chutes \& Ladders in 1943. All variations of the game share the same basic idea; players move between numbered squares by rolling a die (or spinning a spinner) and ascend ladders and slide down chutes throughout the game, all with a goal of getting to the end. Specifically, the official game rules are as follows.
\begin{itemize}
    \item There are 100 Squares on the board. 
    \item The goal is to reach the square-100. 
    \item Players spin a spinner or roll a die to decide how far to move. If you land at the bottom of a ladder you immediately ascend to the top, or if you land at the top of a slide you immediately ride all the way down. Otherwise, you stay where you land. Related pictures connect the end points of each chute and ladder.
    \item Players start off the board. This position may be thought of as square-0.
    \item More than one player may occupy the same square at a time.
    \item To win, a player needs to land on square-100; overshooting it will revert the player to the square they occupied at the start of their turn.
\end{itemize}
Since the player’s motion is random, the entire game is fundamentally a random process. We will construct a Markov Model of the game and run many millions of simulated games to answer questions like,
\begin{itemize}
    \item What is the average duration of the game?
    \item What happens to the average duration of the game as the probability of rolling a $\delta \in \lbrace 1, 2, 3, 4, 5, 6 \rbrace$ approaches 100\%?
    \item What happens to the average duration of the game if we change the rules slightly and introduce strategy to the game?
\end{itemize}

\section{Modelling The Game As A Markov Chain}
\label{Modeling}
Tsitsiklis and  Jaillet provide a phenomenal introduction to Markov Chains in their MIT open course \cite[Lectures 24--25]{MIT_Prob_2018}. Now, in order to model the game as a Markov Chain we will represent each square of the game as a state in the chain. Square-$n$ will be represented by state $n$ at first; we will adjust as needed. To begin, let us consider what is possible at the start of the game; when the player is on square-0. The figure below shows the board where ladders are green and chutes are red.
\begin{figure}[H]
    \centering
    \caption{Chutes \& Ladders Board}
    \includegraphics[scale=0.65]{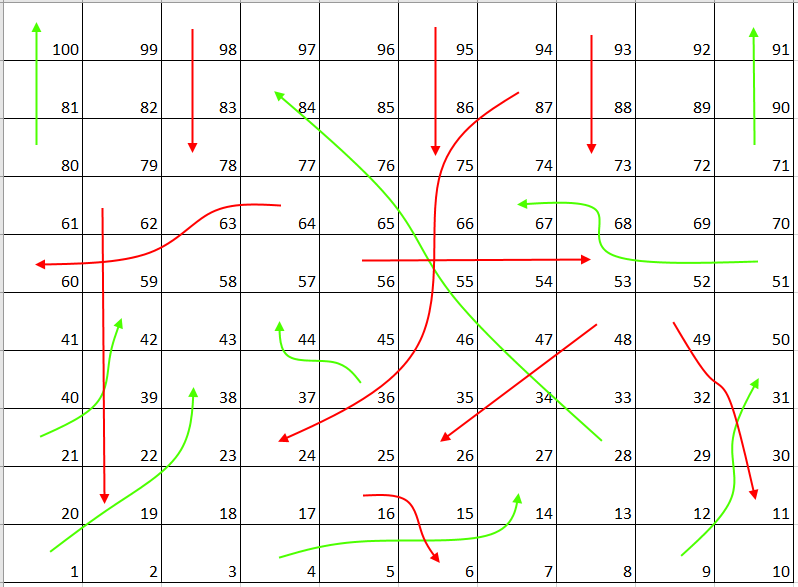}
    \label{board}
\end{figure}
\noindent The player could roll a 1 to advance to square-1 (and they would then advance to square-38), or a 2 to advance to square-2, ..., or a 4 to advance to square-4 (and they would then advance to square-14), ..., or a 6 to advance to square-6. Let $S = \{1,2,3,4,5,6 \}$ be the sample space of all outcomes from rolling a six-sided die, and for each event $\delta$ in $S$, let $P(\delta)$ denote the probability of rolling $\delta$. For a fair die $P(\delta) = 1/6$ for all events $\delta$. Since it is not possible to roll anything higher than a 6, the Markov Chain modeling the game, at state 0, looks something like the following.
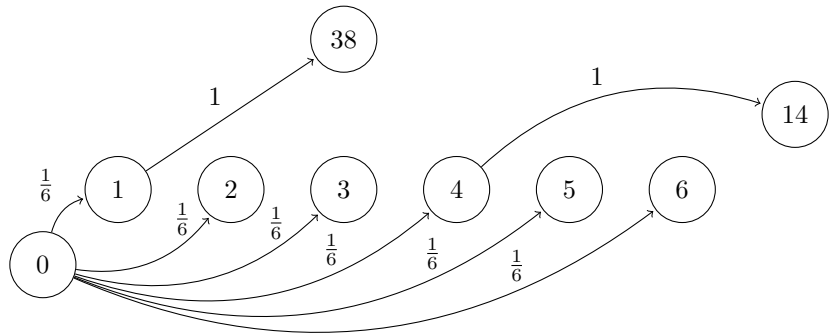
\begin{figure}[H]
    \caption{Markov Model at State 0, Version I}
    \begin{center}
        \begin{tikzpicture}
        \node[state] at (-1,-1)       (0) {0};
        \node[state] at (0,0)        (1) {1};
        \node[state] at (1.5,0)        (2) {2};
        \node[state] at (3,0)        (3) {3};
        \node[state] at (4.5,0)        (4) {4};
        \node[state] at (6,0)        (5) {5};
        \node[state] at (7.5,0)        (6) {6};
        \node[state] at (3,2)        (38) {38};
        \node[state] at (9,1)        (14) {14};
    
    \draw[every loop]
        (0) edge[bend left, auto = left] node {$\frac{1}{6}$} (1)
        (0) edge[bend right, above, pos=0.75] node {$\frac{1}{6}$} (2)
        (0) edge[bend right, above, pos=0.80] node {$\frac{1}{6}$} (3)
        (0) edge[bend right, above, pos=0.70] node {$\frac{1}{6}$} (4)
        (0) edge[bend right, above, pos=0.75] node {$\frac{1}{6}$} (5)
        (0) edge[bend right, above, pos=0.75] node {$\frac{1}{6}$} (6)
        (1) edge[auto = left] node {$1$} (38)
        (4) edge[bend left, auto = left] node {$1$} (14);
        \end{tikzpicture}
    \end{center}
    \label{model1}
\end{figure}
\noindent However, this is not entirely correct. As depicted, if the player were to roll a 1, then on their next turn they are guaranteed to ascend the ladder to square-38. The rules of the game state when the player lands at the base of a ladder they immediately ascend the ladder, and their turn ends thereafter. The same is true regarding chutes. So, the Markov model of the game at state 0 should be:

\begin{figure}[H]
    \caption{Markov Model at State 0, Version II}
    \begin{center}
        \begin{tikzpicture}
            \node[state] at (-1,-1)       (0) {0};
            \node[state] at (1.5,0)        (2) {2};
            \node[state] at (3.5,0)        (3) {3};
            \node[state] at (6,0)        (5) {5};
            \node[state] at (4,2)        (38) {38};
            \node[state] at (9,1)        (14) {14};
            \node[state] at (8,0)        (6) {6};
        
        \draw[every loop]
            (0) edge[bend left, auto = left] node {$\frac{1}{6}$} (38)
            (0) edge[bend right, above, pos=0.7] node {$\frac{1}{6}$} (2)
            (0) edge[bend right, above, pos=0.75] node {$\frac{1}{6}$} (3)
            (0) edge[bend right, above, pos=0.75] node {$\frac{1}{6}$} (6)
            (0) edge[bend right, above, pos=0.60] node {$\frac{1}{6}$} (5)
            (0) edge[out=290,in=270,above,pos=0.5] node {$\frac{1}{6}$} (14);
        \end{tikzpicture}
    \end{center}
    \label{model2}
\end{figure}
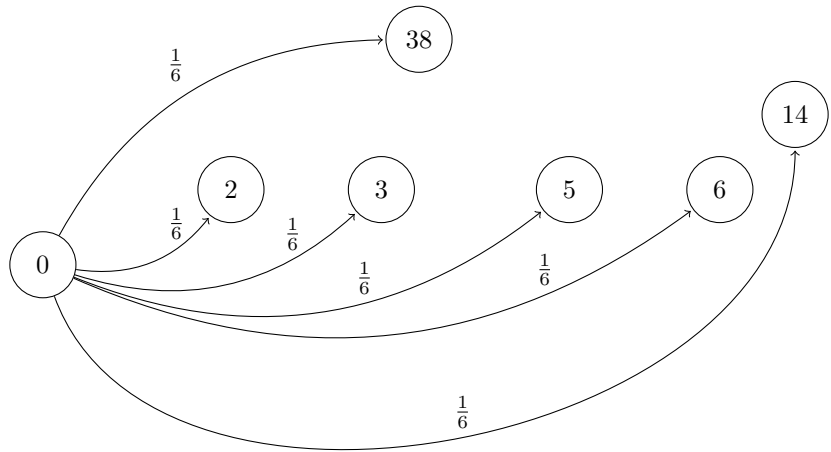

\noindent Notice, it is impossible for the player to roll the die while occupying any square at the base of a ladder or top of a chute. So, no state exists for those squares. Next we need to compile all the transition probabilities into a matrix where the $(i,j)$ entry is the probability of transitioning from state $i$ to state $j$ in one turn. So, from figure \ref{model2},

\makebox[\textwidth][c]{
\begin{tabular}{|c|c|c|c|c|c|c|c|c|c|c|c|c|c|}\hline
    From $\backslash$ To & State 0 & State 2 & State 3 & State 5 & State 6 & State 7 & State 8 & ... & State 14 & ... & State 38 & ... & State 100 \\\hline
    State 0 & 0 & $P(2)$ & $P(3)$ & $P(5)$ & $P(6)$ & 0 & 0 & ... & $P(4)$ & ... & $P(1)$ & ... & 0\\\hline
\end{tabular}
}

\newpage

\begin{figure}[H]
    \caption{Markov Model of First Two Turns, Rolling From 14 and 38 Excluded}
    \begin{center}
    \hspace*{-0.5cm}
        \begin{tikzpicture}
            \node[state] at (-1,-1)       (0) {0};
            \node[state, color=blue] at (1.5,0)        (2) {2};
            \node[state, color=brown] at (3.5,0)        (3) {3};
            \node[state, color=red] at (6,0)        (5) {5};
            \node[state, color=orange] at (4,6)        (38) {38};
            \node[state] at (16,6)        (31) {31};
            \node[state, color=violet] at (9,3)        (14) {14};
            \node[state, color=60color] at (8,0)        (6) {6};
            \node[state] at (10,0)        (7) {7};
            \node[state] at (12,0)        (8) {8};
            \node[state] at (16,0)        (10) {10};
            \node[state] at (16,3)        (11) {11};
            \node[state] at (14,3)        (12) {12};
        
        \draw[every loop]
            (0) edge[bend left, auto = left] node {$P(1)$} (38) %Turn 1
            (0) edge[bend right, above, pos=0.5] node {$P(2)$} (2) %Turn 1
            (0) edge[bend right, above, pos=0.75] node {$P(3)$} (3) %Turn 1
            (0) edge[bend right, below, pos=0.5] node {$P(6)$} (6) %Turn 1
            (0) edge[bend right, above, pos=0.60] node {$P(5)$} (5) %Turn 1
            (0) edge[out=290,in=270,above,pos=0.4, looseness=1.5] node {$P(4)$} (14) %Turn 1
            
            (2) edge[above, pos=0.5, color=blue] node {$P(1)$} (3) %Turn 2
            (2) edge[out = 80, in=120, below, pos=0.75, looseness = 0.6, color=blue] node {$P(3)$} (5) %Turn 2
            (2) edge[out = 120, in=180, above, pos=0.75, looseness = 1, color=blue] node {$P(2)$} (14) %Turn 2
            (2) edge[out = 90, in=150, above, pos=0.75, looseness = 0.8, color=blue] node {$P(4)$} (6) %Turn 2
            (2) edge[out = 100, in=150, above, pos=0.75, looseness = 0.8, color=blue] node {$P(5)$} (7) %Turn 2
            (2) edge[out = 110, in=150, above, pos=0.85, looseness = 0.8, color=blue] node {$P(6)$} (8) %Turn 2

            (3) edge[below, pos=0.5, color=brown] node {$P(2)$} (5) %Turn 2
            (3) edge[out = 60, in=225, above, pos=0.7, looseness = 0.4, color=brown] node {$P(1)$} (14) %Turn 2
            (3) edge[out = 310, in=260, below, pos=0.7, looseness = 0.4, color=brown] node {$P(3)$} (6) %Turn 2
            (3) edge[out = 285, in=235, below, pos=0.8, looseness = 0.7, color=brown] node {$P(4)$} (7) %Turn 2
            (3) edge[out = 280, in=235, below, pos=0.8, looseness = 0.7, color=brown] node {$P(5)$} (8) %Turn 2
            (3) edge[out = 275, in=270, above, pos=0.35, looseness = 1.25, color=brown] node {$P(6)$} (31) %Turn 2

            (5) edge[above, pos=0.3, color=red] node {$P(1)$} (6) %Turn 2
            (5) edge[out=70, in=155, above, pos=0.35, looseness = 0.4, color=red] node {$P(2)$} (7) %Turn 2
            (5) edge[out=275, in=270, above, pos=0.35, looseness = 1.4, color=red] node {$P(3)$} (8) %Turn 2
            (5) edge[out=270, in=260, below, pos=0.4, looseness = 1.8, color=red] node {$P(4)$} (31) %Turn 2
            (5) edge[out=265, in=275, above, pos=0.6, looseness = 1.8, color=red] node {$P(5)$} (10) %Turn 2
            (5) edge[out=260, in=270, above, pos=0.65, looseness = 1.5, color=red] node {$P(6)$} (11) %Turn 2

            (6) edge[below, pos=0.5, color=60color] node {$P(1)$} (7) %Turn 2
            (6) edge[out=60, in=155, above, pos=0.4, looseness = 0.4, color=60color] node {$P(2)$} (8) %Turn 2
            (6) edge[out=105, in=200, above, pos=0.5, looseness = 0.4, color=60color] node {$P(3)$} (31) %Turn 2
            (6) edge[out=100, in=120, above, pos=0.60, looseness = 0.4, color=60color] node {$P(5)$} (11) %Turn 2
            (6) edge[out=90, in=125, above, pos=0.6, looseness = 0.55, color=60color] node {$P(4)$} (10) %Turn 2
            (6) edge[out=95, in=210, above, pos=0.60, looseness = 0.4, color=60color] node {$P(6)$} (12); %Turn 2
        \end{tikzpicture}
    \end{center}
    \label{model3}
\end{figure}
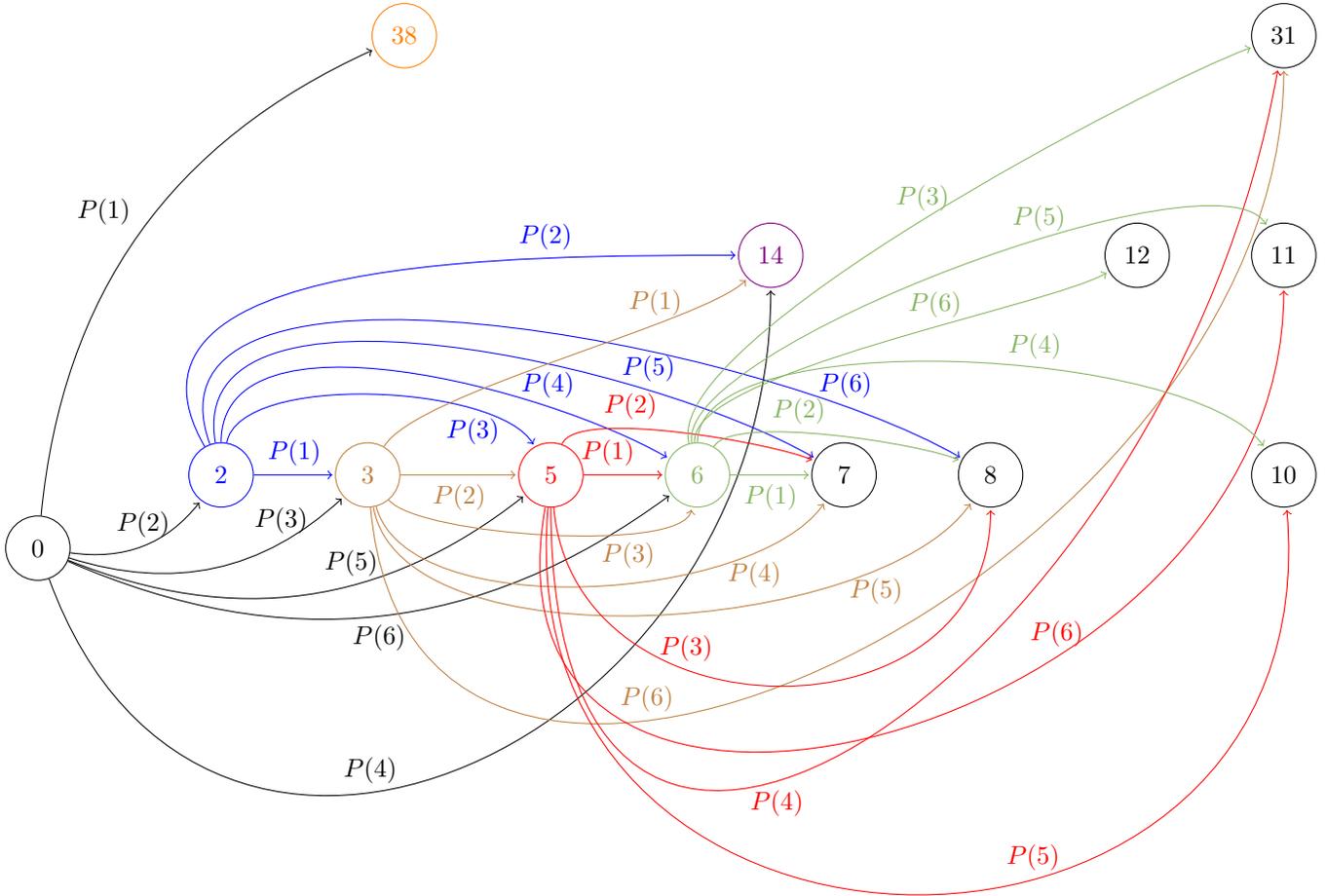

Figure \ref{model3} shows the Markov Chain that models the player's first two turns of the game.

\noindent All possible transitions for the Player's firs turn are shown in black. For the sake of simplicity, we omit from Figure \ref{model3} the possible transitions on the Player's second turn from square-$38$ and from square-$14$. For any two states, $i$ and $j$, in figure \ref{model3} such that there does not exist a transition $i \to j$, it is implied $P(i \to j)$ is zero. Thus, each state $i$, corresponding to row $i$ of the transition probability table, will have exactly six nonzero entries. That is, row $i$ will have the value $P(j)$ in column $i+j$ unless square-$(i+j)$ is at the bottom of a ladder or the top of a chute. In that case, row $i$ will have value $P(j)$ in the column corresponding to the termination square of the chute or ladder beginning at square-$(i+j)$. Then, since square-1 is at the bottom of a ladder it is not possible to roll from square-1. So, the next row in the table is

\makebox[\textwidth][c]{
\begin{tabular}{|c|c|c|c|c|c|c|c|c|c|c|c|c|c|}\hline
    From $\backslash$ To & State 0 & State 2 & State 3 & State 5 & State 6 & State 7 & State 8 & ... & State 14 & ... & State 38 & ... & State 100 \\\hline
    State 0 & 0 & $P(2)$ & $P(3)$ & $P(5)$ & $P(6)$ & 0 & 0 & ... & $P(4)$ & ... & $P(1)$ & ... & 0\\\hline
    State 2 & 0 & 0 & $P(1)$ & $P(3)$ & $P(4)$ & $P(5)$ & $P(6)$ & ... & $P(2)$ & ... & 0 & ... & 0\\\hline
\end{tabular}
}

\noindent The row State 2 corresponds to all transitions in figure \ref{model3} of the form $2 \to j$. After doing this for all 82 possible starting states we will have constructed quite a sparse matrix $P$. The final row of which is,

\makebox[\textwidth][c]{
\begin{tabular}{|c|c|c|c|c|c|c|c|c|c|c|c|c|c|}\hline
    From $\backslash$ To & State 0 & State 2 & State 3 & State 5 & State 6 & State 7 & State 8 & ... & State 14 & ... & State 38 & ... & State 100 \\\hline
    State 100 & 0 & 0 & 0 & 0 & 0 & 0 & 0 & ... & 0 & ... & 0 & ... & 1\\\hline
\end{tabular}
}

\noindent Since when the player reaches square-100 the game is over and so they stay on square-100. With this matrix we may begin answering some questions about the game.
\subsection{Calculating the Probability the Game Will Last $M$ Rounds} 
\begin{figure}[H]
    \caption{$P(M)$ For $1 \leq M \leq 100$}
    \label{P(M)}
    \begin{center}
        \includegraphics[scale=0.75]{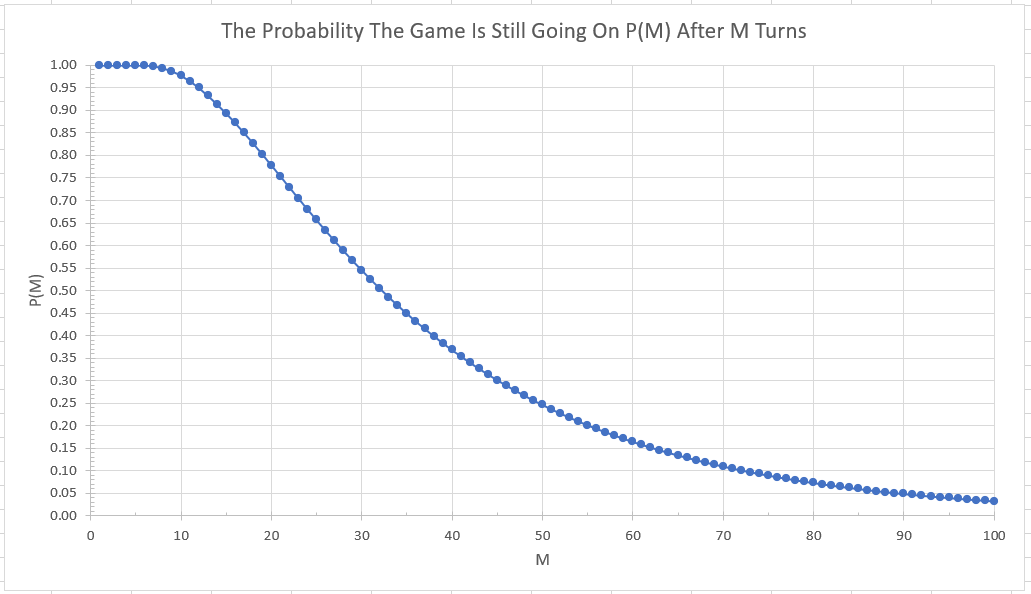}
\end{center}
\end{figure}

To calculate the probability the game will still be going on after $M$ rounds we restrict our attention to $T = P_{transient}$. That is, model only the transient states of the game; or remove all rows and columns associated with any absorbing state. For the case of Chutes \& Ladders, this means removing the 82nd column and the 82nd row; as square-100 is the only absorbing state. Matrix $T$ is shown in Listing 2 (figure \ref{sm1}). The matrix is far to large to view properly in this PDF; we suggest the reader copy \& paste the code in Listing 2 into the IDE of their choice. We recommend Cocalc SageMath. Next we must use the following equation
\[P(M) = e_s \cdot T^M \cdot \mathbbm{1}.\]
Where $e_s \in \mathbbm{R}^{81}$ represents the players starting position by having a $1$ in position $s$ and $0$ elsewhere; note, $s$ need not represent square-$s$. For example $s=1$ is square-0, and $s=81$ is square-99. Moreover, $\mathbbm{1} \in \mathbbm{R}^{81}$ such that every entry is a $1$. For example, starting from square-0 the probability the game is still going on after $M$ turns is shown in figure \ref{P(M)}.

\subsection{Calculating The Expected Duration of The Game $\sigma$}\label{avg} The probability the game is still going on after $M$ rounds goes to zero as $M$ goes to infinity. Since each state of matrix $T$ is occupied by the player only one turn at a time the expected duration of the game is the sum of the probability the game is still going on, over all possible number of turns. That is,
\begin{align*}
    \text{Avg. Duration} &= \sum_{M=0}^\infty e_s \cdot T^M \cdot \mathbbm{1}\\
    &= e_s \cdot \left[\sum_{M=0}^\infty T^M\right] \cdot \mathbbm{1}\\
    &= e_s \cdot \left[I + T + T^2 + T^3 + \dots\right] \cdot \mathbbm{1}\\
    &= e_s \cdot \left[ I-T \right]^{-1} \cdot \mathbbm{1}
\end{align*}
Where $I$ is the identity matrix. We used the above equation in another paper involving Markov Models \cite{bardenova2021markov}. Using a fair die and starting from square-0 ($e_1$), we calculate,
\[ e_1 \cdot (I-T)^{-1} \cdot \mathbbm{1} = 39.5984 \]
That is, the average duration of a standard game of Chutes \& Ladders is about 39.5 turns. If each turn takes approximately thirty seconds, then the game would last about twenty minutes.

\section{Fully Weighted Die Cases}
\label{Full}
There are some mathematical curiosities to be noted when considering the case where the player is using a die that has nearly 100\% change of rolling one of 1, 2, 3, 4, 5, or 6. Yet first, we need some definitions.
\begin{enumerate}
    \item Let the average number of turns per game be denoted as $\sigma$.
    \item Let a $\delta$\textbf{-chute-loop} or $\delta$\textbf{-c-loop} be a set that contains all squares such that, for $P(\delta) = 1$, should the player land on a square in a $\delta$-c-loop, they will indefinitely move between all squares in the $\delta$-c-loop on a fixed repeating path. An example of a 5-c-loop is $\lbrace 6, 11, 16 \rbrace$.
    \item Let a $\delta$\textbf{-finale-loop} or $\delta$\textbf{-f-loop} be a set containing a single square such that, for $P(\delta) = 1$, should the player land on this square a roll of $\delta$ is too large to win the game. An example of a 4-f-loop would be $\lbrace 97 \rbrace$. Note: this set is also a 5-f-loop and 6-f-loop.
    \item Let $\mathbf{\delta}$\textbf{-loop} refer to either a $\delta$-c-loop or a $\delta$-f-loop.
    \item Let the $\mathbf{\delta}$\textbf{-win path} be the set of all squares such that, for $P(\delta) = 1$, if the player were to land on any square within the $\delta$-win path they would be guaranteed to win the game in a finite number of moves. For example, the $6$-win path is $\lbrace 68, 70, 74, 76, 80, 82, 88, 94 \rbrace$. Notice $80 \in \delta$-win path for all $\delta$ since square-80 is at the base of a ladder that immediately takes you to square-100.
    \item Let a $\delta$\textbf{-on ramp} be the set of squares such that, for $P(\delta = 1)$, should the player land on any square in a $\delta$-on ramp they are guaranteed to eventually either enter a $\delta$-c-loop, or a $\delta$-f-loop. Different $\delta$-loops have different $\delta$-on ramps. An example of a 5-on ramp is $\lbrace93, 98\rbrace$. That 5-on ramp is guaranteed to lead the player to the 5-c-loop $\lbrace 73, 78, 83, 88 \rbrace$.
\end{enumerate}

\subsection{Fully Weighted on 1,2,3, or 6}
\label{on1236}
\begin{figure}[H]
    \centering
    \caption{Square Classification for $P(3)=1$}
    \includegraphics[scale=0.77]{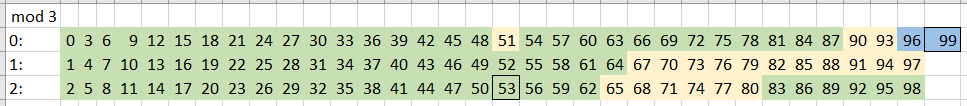}
    \label{3class}
\end{figure}

If the player uses a die that has a 100\% chance of rolling one of 1, 2, 3, or 6, they will always get stuck in a $\delta$-loop. In fact, for $\delta \in \lbrace 1,2,3,6\rbrace$, as the probability of rolling a $\delta$, $P(\delta)$, approaches 100\% the average number of turns per game diverges. $$\ds \lim_{P(\delta) \to 1}\sigma = \infty$$ Yet, what may not be obvious is that as $P(3) \to 1$, $\sigma$ ascends to infinity much faster than for any other $\delta \in \lbrace 1,2,6 \rbrace$. Below is a graph of the average number of turns per game vs the probability of rolling each number $\left(P(\delta)\right)$. Each average was calculated by simulating one hundred thousand games. Thus, each point in the graph below represents the average number of turns over 100,000 games. These are known as 'Monte Carlo' simulations, the principles and history of which are reviewed by Harrison \cite{Harrison2010}.
\begin{center}
    \includegraphics[scale=0.6]{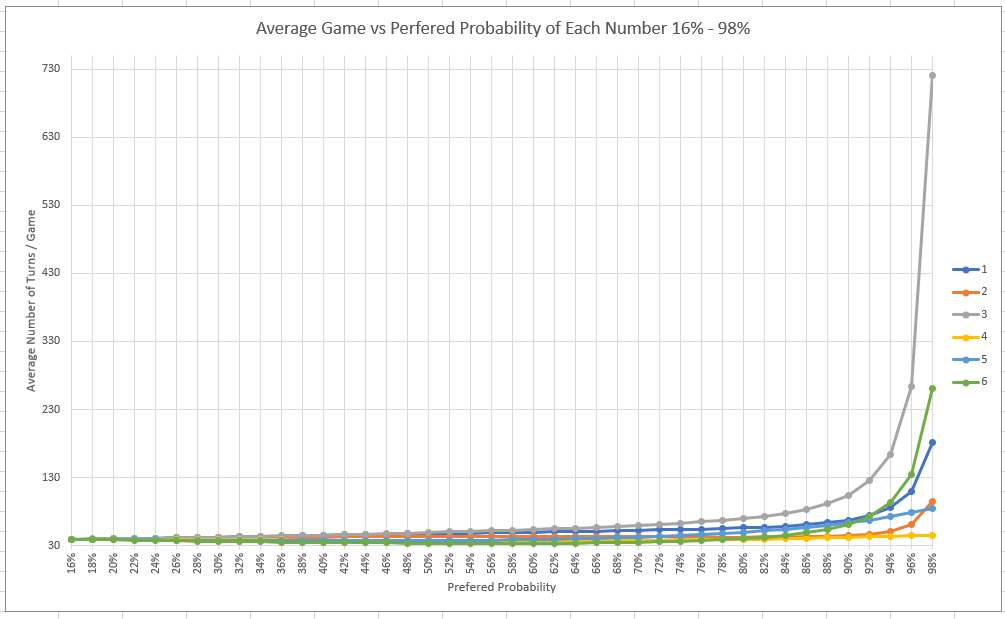}
\end{center}
Clearly, as $P(3) \to 1$, $\sigma \to \infty$ much faster than the rest. There is at least one major factor contributing to this; rolling mostly all 3s gets the player stuck in a loop most of the time. In fact, letting the player start from each square and noting what happens if $P(3) = 1$, it may be calculated that if the player were to select a random square on the board to start from then 80\% of the time they will be stuck in an  3-loop. This is the highest probability of all $\delta$-loops. That is, if the player were to select a random square on the board to start from, where $P(\delta) = 1$, the probability that the player will be in a $\delta$-loop is highest at $\delta = 3$. Indeed, among all 100 possible starting squares (0 - 99) the number of squares that fall into each classification, $\delta$-C-Loop, $\delta$-F-Loop, or $\delta$-Win Path, all including their respective $\delta$-on ramps, are as follows.
\begin{center}
    \begin{tabular}{|c|c|c|c|}\hline
         \multicolumn{4}{|c|}{$\delta$ vs Number of Squares per Classification} \\\hline
         $P(\delta)=1$ & $\delta$-C-Loop & $\delta$-F-Loop & $\delta$-Win Path\\\hline
         $\delta=1$ & 70 & 0 & 30\\\hline
         $\delta=2$ & 9 & 2 & 89\\\hline
         $\delta=3$ & 78 & 2 & 20\\\hline
         $\delta=4$ & 34 & 2 & 74\\\hline
         $\delta=5$ & 54 & 23 & 23\\\hline
         $\delta=6$ & 0 & 92 & 8\\\hline
    \end{tabular}
\end{center}
Notably, for $P(1)=1$, if the player were to select a random square on the board to start from, they will never be trapped in a 1-finale loop. This is of course because it is impossible to overshoot square-100 if the player always rolls a 1. Surprisingly, no matter where the player starts on the board, if $P(6)=1$, they will never get stuck in a 6-c-loop, yet 92\% of the time they'll be stuck in an 6-f-loop. Intuitively, it may seem that for $P(\delta)=1$, since $\delta = 6$ has a higher chance of being in a 6-loop than $\delta = 3$ has of being in a 3-loop, as $P(6) \to 1$, $\sigma \to \infty$ faster than $\sigma \to \infty$ as $P(3) \to 1$. Yet, this is not the case. For $\delta = 3$, there actually exists only one 3-c-loop; for which, 77 different squares makeup its 3-on ramp. That 3-c-loop contains solely square-53. This will be shown in Figure $\ref{3class}$ in Section \ref{classes}. The only way the player may escape square-53 involves rolling more than one non-3 in quick succession; which is extremely unlikely as $P(3) \to 1$. However, when $\delta = 6$, the player is only ever stuck in a 6-f-loop. So, they only ever need to roll one non-6 to win the game when stuck in a 6-f-loop. This is clearly much more likely as $P(6) \to 1$ then rolling two non-3s in quick succession is as $P(3) \to 1$. Furthermore, using the markov model of the game shown in Section \ref{avg} we have calculated that for $P(6) = 0.99999999999999999999999$ (all but one in one hundred sextillion), $\sigma \approx 5.0E+23$ where as for $P(3)$ of the same probability $\sigma \approx 1.9E+45$. Thus, certainly as $P(3) \to 1$, $\sigma \to \infty$ much faster than for any other $\delta$.

Another interesting item is, of all $\delta$ such that $\sigma \to \infty$ as $P(\delta) \to 1$, $\sigma$ appears to ascend to infinity the slowest as $P(2) \to 1$. This is obvious from the graph below which is the same graph as above, just zoomed in to see detail. Note: for $\delta = 4,5$, $\sigma \not \to \infty$ as $P(\delta)\to 1$ (we investigate these in great detail later).

\begin{center}
    \includegraphics[scale = 0.6]{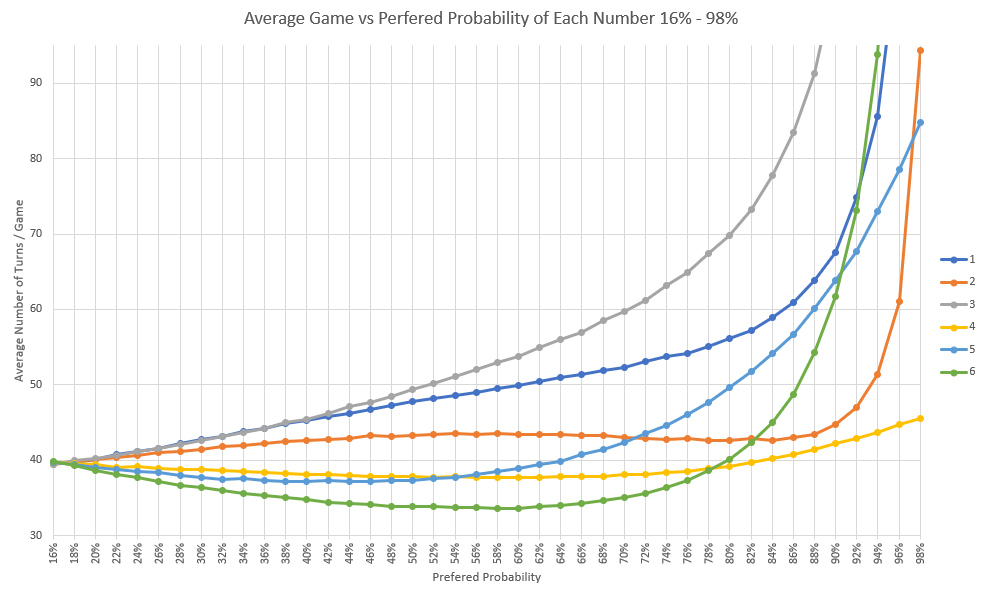}
\end{center}

Lastly, it appears that for the shortest game, the player should use a die that is weighted to roll a 6 about 60\% of the time. In the next section we consider the classification of each square on the board, for $P(\delta) = 1$, for every $\delta$. That is, should the player land on square-$x$ with $P(\delta)=1$, what happens? Does the player get stuck in a $\delta$-loop or do they win?

\subsection{Square Classifications}\label{classes} Below are the classifications of every square on the board, given $P(\delta)=1$. Squares shaded the same color are contained inside a certain $\delta$-loop or its $\delta$-on ramp. If a square is inside of an outlined area then that square is in the $\delta$-loop for that color, and if a square is not inside an outlined are then that square is in the $\delta$-on ramp for that color. The $\delta$-winning path is always colored \includegraphics[scale=0.25]{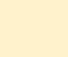}.

\begin{figure}[H]
    \centering
    \caption{Square Classification for $P(1)=1$}
    \includegraphics{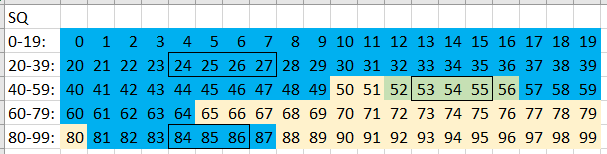}
    \label{1class}
\end{figure}
\noindent So, from figure \ref{1class}, square-52 and square-56 are contained inside the $1$-on ramp for the $1$-c-loop $\lbrace 53,54,55 \rbrace$. That is, should the player land on square-52 or square-56 with $P(1) = 1$, then they will get stuck in an infinite loop traveling between squares $\lbrace 53,54,55 \rbrace$. Recall, Figure \ref{board} to see why this is true.

\begin{figure}[H]
    \centering
    \caption{Square Classification for $P(2)=1$}
    \includegraphics{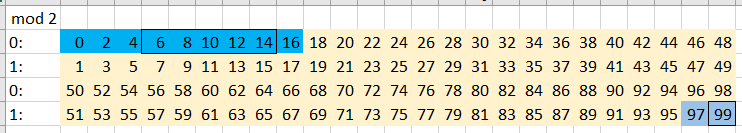}
    \label{2class}
\end{figure}
\noindent Notice, for $P(2) = 1$ the board has the most number of squares in the $\delta$-win path of all $\delta$.

\noindent For $P(3)=1$, there is only once $3$-c loop and one $3$-f loop; $\lbrace 53 \rbrace$ and $\lbrace 99 \rbrace$ respectively.

\begin{figure}[H]
    \centering
    \caption{Square Classification for $P(6)=1$}
    \includegraphics{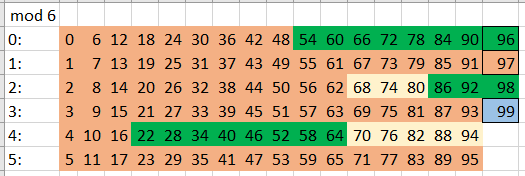}
    \label{6class}
\end{figure}
\noindent Of all $\delta$, $P(\delta=6)=1$ is, the only one that does not have any $\delta$-c loops.

\subsection{Fully Weighted on 5}\label{on5} For $P(\delta = 5) = 1$, the player will win the game in 16 turns every time. That is, if $P(5)=1$ then $\sigma = 16$. Yet, for P(5) very close to 1, and $P(\delta \not = 5) = [1-P(5)]/5$, it may be calculated using the Markov model of the game that
\[\sigma = e_1 \cdot (I - T)^{-1} \cdot \mathbbm{1} \approx 82.37\]
So, it appears that $\lim_{P(5) \to 1} \sigma \approx 82.37$. Why is this the case? Why doesn't $\sigma \to 16$ as $P(5) \to 1$? The average number of turns per game may be approximated using figure \ref{5class} below. Squares contained in the $\delta$-win path that the player will transition between if they start on square-0 are outlined. It will be necessary to distinguish between these and other squares contained in the $\delta$-win path to model $\lim_{P(\delta) \to 1} \sigma$. 

\begin{figure}[H]
    \centering
    \caption{Square Classification for $P(5)=1$}
    \includegraphics{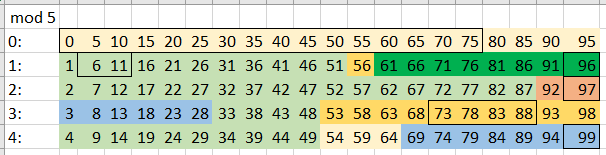}
    \label{5class}
\end{figure}

Consider the infinite game, that is the game such that when the player reaches square-100 they immediately return to square-0 and play again. From the figure above, if $P(5) = 1$ the player will always follow the outlined $5$-win path. Yet, if $P(5)$ were very close to 1, the player would spend most turns in the outlined 5-win path, but not all turns. So, how would the player stay in the 5-win path, if they are already somewhere in the outlined 5-win path? For square-0 the player stays in the 5-win path only if they roll 5, any other roll will cause them to transition to a 5-loop, or the 5-on ramp of a 5-loop. This is also the case for square-5. This is the case for every square in the outlined 5-win path except square-50, square-55, and square-60; for these squares the player will say in the 5-win path if they roll a 5 or a 4. Note, since $P(5)$ is very close to 1, most of the time, if the player enters a $5$-on ramp, they will soon thereafter enter its respective 5-loop. Also, in the infinite game, the outlined $\delta$-win path acts similar to a $\delta$-loop. So, it's a fair assumption that if the player rolls a 4 from square-50 then they will return to to the outlined 5-win path (we are assuming since $P(5)$ is very close to one, the player never sees two or more non-5 rolls in quick succession). If the player is somewhere in the outlined 5-win path then the probability the player is on square-50, square-55, or square-60 is $\frac{3}{16}$. Furthermore, if the player does not roll a 5 when on one of those squares the probability of staying on the 5-win path is $\frac{1}{5}$ since only one non-5 will keep them in the 5-win path (rolling a 4). So, the total probability the player transitions from the outlined 5-win path to the outlined 5-win path in one turn is 
$$P(Win \to Win) = P(5) + \frac{3}{16}\cdot\frac{1}{5}\cdot[1-P(5)] = P(5) + \frac{3}{80}\cdot[1-P(5)]$$
That is, they either need to roll a $5$ or $3/16$ of the time they are also permitted to roll a $4$, which represents $1/5$ of all non-5 rolls. Then, since the player must either be in the outlined 5-win path or a 5-loop, the probability of transitioning from the outlined 5-win path to a 5-loop is
$$P(Win \to Loop) = 1-P(Win \to Win) = 1 - \left(P(5) + \frac{3}{80}\cdot[1-P(5)]\right) = \frac{77}{80}\cdot[1-P(5)]$$
Conversely, if the player is in any $5$-loop, if they roll a $5$ then, by definition, they stay in that 5-loop. However, from the figure above, for all squares in all 5-loops, there is only one non-5 roll that will allow the player to escape back to the 5-win path. So, the probability of transitioning from a 5-loop to the outlined 5-win path is
$$P(Loop \to Win) = \frac{1}{5}\cdot[1-P(5)]$$
Then, the probability of transitioning from a 5-loop to a 5-loop is
$$P(Loop \to Loop) = 1 - P(Loop \to Win) =  P(5) + \frac{4}{5}\cdot[1-P(5)]$$
With these probabilities we may build the Markov model shown in figure \ref{tp_5} below.

\begin{figure}[H]
    \caption{Markov Chain for the Player’s Transition Between the 5-Win Path and 5-loops}
    \label{tp_5}
    \begin{center}
        \begin{tikzpicture}
                % Add the states
                \node[state] at (-3,0)       (w) {Win};
                \node[state] at (3,0)       (l) {Loop};

                \draw[every loop]
                    (w) edge[bend left, auto=left] node {$\frac{77}{80} \cdot [1-P(5)]$}(l)
                    (w) edge[loop left] node {$P(5) + \frac{3}{80} \cdot [1-P(5)]$}(w)
                    (l) edge[bend left, auto=left] node {$\frac{1}{5} \cdot [1-P(5)]$}(w)
                    (l) edge[loop right] node {$P(5) + \frac{4}{5} \cdot [1-P(5)]$}(l);
             \end{tikzpicture}
    \end{center}
\end{figure}

\noindent Note, this markov chain does not model the game in its entirety as it does not consider the possibility of non-5 rolls while on any on-ramp because that would require non-5 rolls in quick succession, which we are ignoring to simplify the model. In reality, it is possible for a player to roll a non-5, placing them on some 5-on ramp and before reaching that 5-on ramp's 5-loop, roll another non-5. However, as it may be obvious, this become exceedingly rare as $P(5) \to 1$. Thus, the model approximates the game well for $P(5)$ close to 1. 

Below we calculate the proportion of turns during the infinite game the player spends on both the 5-win path and a 5-loop using this model. Let $p_w$ and $p_l$ represent the probability that after many moves the player may be found in the "Win" state and in the "Loop" state respectively. This is equivalent to the proportion of turns spent in either state. Then,
\begin{align}
    p_w &= p_wP(Win \to Win) + p_lP(Loop \to Win)\\
        &= p_w \left[ P(5) + \frac{3}{80}\cdot\left[1-P(5)\right] \right] + p_l \left[ \frac{1}{5}\cdot \left[1-P(5)\right] \right]\\
    p_l &= p_wP(Win \to Loop) + p_lP(Loop \to Loop)\\
        &= p_w\left[ \frac{77}{80}\cdot\left[1-P(5)\right] \right] + p_l\left[ P(5) + \frac{4}{5}\cdot \left[ 1-P(5)\right] \right]
\end{align}
In words, the probability the player is on the 5-win path $(p_w)$ is equal to the probability they were already there and then they stay there after the next turn, plus the probability they were on a 5-loop and transition to the 5-win path after the next turn. The description of $p_l$ is similar. Finally, one last equation is needed. For every turn the player must be somewhere. Therefore, it must be such that
\begin{align}
    1 &= p_w + p_l
\end{align}
Now, the system of three equations with two unknowns may be solved. With $P(5) = p(x) =\frac{x-1}{x}$ and $1-P(5) = 1-p(x) = \frac{1}{x}$,
\begin{align}
    p_w &= 1- p_l & \text{From (5)}\\
    p_w &= p_w\left(\frac{x-1}{x} + \frac{3}{80}\cdot\frac{1}{x}\right) + p_l\left( \frac{1}{5}\cdot\frac{1}{x} \right) & \text{From (2)}\\
    &=\left(1-p_l\right)\left(\frac{80x-77}{80x}\right) + \frac{1}{5x}p_l & \text{From (6) and (7)}\\
    &= \frac{80x-77}{80x} - \frac{80x-77}{80x}p_l + \frac{16}{80x}p_l & \\
    &=\frac{80x-77}{80x} + \frac{77-80x}{80x}p_l + \frac{16}{80x}p_l & \\
    &= \frac{80x-77}{80x} + \frac{93-80x}{80x}p_l &\\
    p_l &= p_w\left( \frac{77}{80}\cdot\frac{1}{x} \right) + p_l\left( \frac{x-1}{x} + \frac{4}{5}\cdot\frac{1}{x} \right) & \text{From (4)}\\
    &= p_w\left(\frac{77}{80x}\right) + p_l\left(\frac{5x-1}{5x}\right)\\
    &= \left( \frac{80x-77}{80x} + \frac{93-80x}{80x}p_l \right)\left(\frac{77}{80x}\right) + p_l\left(\frac{5x-1}{5x}\right) & \text{From (11)}\\
    &= \frac{6160x-5929}{80^2x^2} + \frac{7161 - 6160x}{80^2x^2}p_l + \frac{80^2x^2-1280x}{80^2x^2}p_l &\\
    &= \frac{6160x-5929}{80^2x^2} + \frac{80^2x^2-7440x+7161}{80^2x^2}p_l &\\
    \implies p_l &- \frac{80^2x^2-7440x+7161}{80^2x^2}p_l = \frac{6160x-5929}{80^2x^2} & \text{Since $p_l = (16)$}\\
    &\iff \left( 1 -  \frac{80^2x^2-7440x+7161}{80^2x^2}\right)p_l = \frac{6160x-5929}{80^2x^2} &\\
    &\iff \left( \frac{80^2x^2}{80^2x^2} -  \frac{80^2x^2-7440x+7161}{80^2x^2}\right)p_l = \frac{6160x-5929}{80^2x^2} &\\
    &\iff \frac{7440x-7161}{80^2x^2}p_l = \frac{6160x-5929}{80^2x^2} &\\
    \iff p_l &= \frac{6160x-5929}{7440x-7161}
\end{align}
It may not be immediately noticeable, yet it is true that
\[ 6160x-5929 = \frac{77}{93}\left( 7440x-7161 \right) \]
Thus, (21) may be rewritten as
\begin{align}
    p_l &= \frac{\frac{77}{93}\left(7440x-7161\right)}{7440x-7161}\\
    \implies p_l &= \frac{77}{93}
\end{align}
And thus,
\begin{align}
    p_w &= 1 - \frac{77}{93} = \frac{16}{93} 
\end{align}
Therefore, during the infinite game the player may be found in a 5-loop about $77/93 \approx 82.8\%$ of the time and in the 5-win path about $16/93 \approx 17.2\%$ of the time. Curiously, these two proportions do not depend on the value of $P(5)$. 

Finally, for a single game and for, $$P(5) = p(x) = \frac{x-1}{x}$$
the probability the game ends in 16 moves is $p(x)^{16}$. This contributes 16 turns to the expected number of turns. The probability this does not happen is $\left[ 1-p(x)^{16} \right]$. In this case, according to our model the player spends about $82.8\%$ of turns in a 5-loop. Since, on average, the player rolls a non-5 only once every $x$ turns and only one non-5 will let them escape back to the 5-win path, when the player is in a 5-loop they will need to roll $5x$ times until they escape to the 5-win path. So, this case contributes $\frac{77}{93}\cdot5x$ turns to the expected number of turns. Hence, the expected number of turns per game is
\begin{align*}
    E_5(x) &= p(x)^{16} \cdot 16 + \left[ 1 - p(x)^{16} \right] \cdot \frac{77}{93} \cdot 5x\\
    &\approx p(x)^{16} \cdot 16 + \left[ 1 - p(x)^{16} \right] \cdot 4.14x
\end{align*}
Where $E_5(x)$ is a function of the average number of rolls per single non-5 roll. That is, the number $x$ in the equation above is the average number of turns per non-5 roll. Notice, the average duration of the game for $P(5)=1$ is approximately,
\[ \sigma \approx \lim_{x \to \infty} E_5(x) \approx 82.24\]
Clearly, $82.24 \not = 82.37$; that is, our approximation of $\sigma$ from $E_5$ does not equal true $\sigma$. This is because the model is not perfect; it does not account for the possibility of more than one non-5 roll in quick succession. However, the model is within $0.16\%$ of true $\sigma$. This model demonstrates why $P(5) \to 1 \implies \sigma \to 82.37$ whereas $P(5)=1 \implies \sigma = 16$. For $P(5)$ close to 1, most games end in $16$ moves. However, over a very large number of games, the player spends about $82.8\%$ of all turns in a 5-loop and only about $17.2\%$ of all turns in the 5-win path. So, any game not ending in 16 moves, on average, will take an extremely large number of turns to complete.

\subsection{Fully Weighted on 4}\label{on4} 
\begin{figure}[H]
    \centering
    \caption{Square Classification for $P(4)=1$}
    \includegraphics{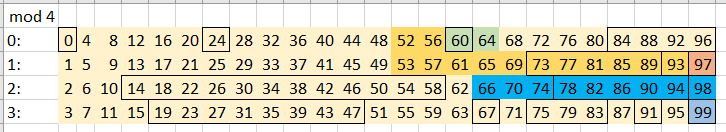}
    \label{4class}
\end{figure}

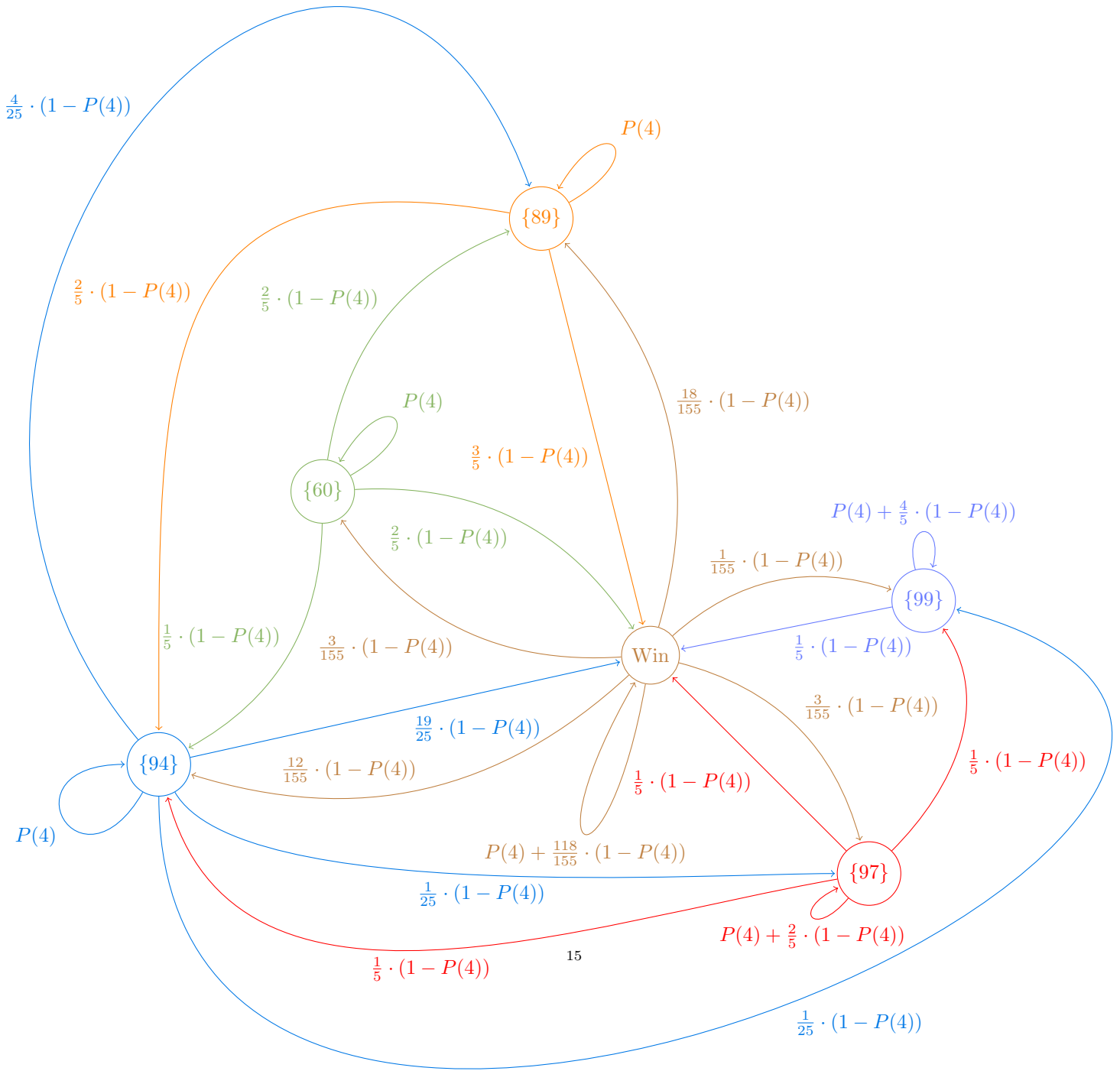
\begin{figure}[H]
    \caption{Markov Chain for the Player's Transition Between the 4-Win Path and 4-loops}
    \label{tp_4}
\begin{changemargin}{-5.25cm}{0cm}

\begin{tikzpicture}[scale=.95]
        % Add the states
        \node[state, brown] at (2,-3)       (w) {Win};
        \node[state, 99color] at (7,-2)       (99) {$\lbrace 99 \rbrace$};
        \node[state, red] at (6,-7)       (97) {$\lbrace 97 \rbrace$};
        \node[state, 60color] at (-4,0)       (60) {$\lbrace 60 \rbrace$};
        \node[state, 94color] at (-7,-5)       (94) {$\lbrace 94 \rbrace$};
        \node[state, orange] at (0,5)       (89) {$\lbrace 89 \rbrace$};

        \draw[every loop]
            (w) edge[bend left, auto=left, color=brown] node {$\frac{3}{155} \cdot (1-P(4))$}(60)
            (w) edge[bend left, auto=left, color=brown] node {$\frac{3}{155} \cdot (1-P(4))$}(97)
            (w) edge[bend left, midway, above, color=brown] node {$\frac{1}{155} \cdot (1-P(4))$}(99)
            (w) edge[bend left, auto=right, color=brown] node {$\frac{12}{155} \cdot (1-P(4))$}(94)
            (w) edge[bend right, auto=right, color=brown] node {$\frac{18}{155} \cdot (1-P(4))$}(89)
            (w) edge[out=260,in=240,loop, looseness = 55, midway, below, color=brown] node {$P(4) + \frac{118}{155} \cdot (1-P(4))$} (w)

            (99) edge[auto=left,color=99color] node {$\frac{1}{5} \cdot (1-P(4))$}(w)
            (99) edge[loop above,color=99color] node {$P(4) + \frac{4}{5} \cdot (1-P(4))$}(99)

            (60) edge[bend left, auto=right, color=60color] node {$\frac{2}{5} \cdot (1-P(4))$}(w)
            (60) edge[bend left, auto=left, color=60color] node {$\frac{2}{5} \cdot (1-P(4))$}(89)
            (60) edge[bend left, auto=right, color=60color] node {$\frac{1}{5} \cdot (1-P(4))$}(94)
            (60) edge[loop, out=30, in=60, looseness = 15, auto=right, color=60color] node {$P(4)$}(60)

            (89) edge[auto=right, color=orange] node {$\frac{3}{5} \cdot (1-P(4))$}(w)
            (89) edge[out=170, in=90, looseness=1.5, auto=right, color=orange] node {$\frac{2}{5} \cdot (1-P(4))$} (94)
            (89) edge[loop, out=30, in=60, looseness = 15, auto=right, color=orange] node {$P(4)$}(89)

            (94) edge[auto=right, color=94color] node {$\frac{19}{25} \cdot (1-P(4))$}(w)
            (94) edge[out=130, in=110, looseness=2, auto=left, color=94color] node {$\frac{4}{25} \cdot (1-P(4))$} (89)
            (94) edge[out=300, in=180, looseness=0.5, auto=right, color=94color] node {$\frac{1}{25} \cdot (1-P(4))$} (97)
            (94) edge[out=270, in=345, looseness=2, auto=right, color=94color] node {$\frac{1}{25} \cdot (1-P(4))$} (99)
            (94) edge[out=240, in=180, loop, auto=left, looseness=10, color=94color] node {$P(4)$}(94)

            (97) edge[auto=left, color=red] node {$\frac{1}{5} \cdot (1-P(4))$} (w)
            (97) edge[out=45,in=305,auto=right, color=red] node {$\frac{1}{5} \cdot (1-P(4))$} (99)
            (97) edge[out=190,in=285,auto=left, color=red] node {$\frac{1}{5} \cdot (1-P(4))$} (94)
            (97) edge[out=230, in=205, loop, midway, below, looseness=10, color=red] node {$P(4) + \frac{2}{5}\cdot(1-P(4))$}(97);
            
     \end{tikzpicture}
\end{changemargin}
\end{figure}

A similar but more complicated case happens for $P(\delta=4)=1$; the player wins the game every time in $31$ moves. Yet, for P(4) very close to 1, and $P(\delta \not = 4) = [1-P(4)]/5$, it may be calculated using the Markov model of the full game that
\[\text{true }\sigma = e_1 \cdot (I - T)^{-1} \cdot \mathbbm{1} \approx 46.98\]
To find an equation for the expected number of turns per game, $E_4(x)$, we may use figure \ref{4class} to build a similar Markov model as we did in section \ref{on5}. However, we need to model each 4-loop as its own state since since the probability of transitioning from a 4-loop to the outlined 4-win path depends on which 4-loop the player is in.
\noindent Using the figure above, we write down all ways all possible transitions may occur. For example, with the last square of each 4-loop acting a representative for that loop and $\lbrace w \rbrace$ representing the 4-win path, the following are all necessary die rolls to transition to all possible states, in one turn, starting from all squares in the 4-c loop $\lbrace94\rbrace$.
\begin{align*}
    &\lbrace 94\rbrace\to\lbrace w \rbrace & &\lbrace 94\rbrace\to\lbrace 60 \rbrace & &\lbrace 94\rbrace\to\lbrace 94 \rbrace & &\lbrace 94\rbrace\to\lbrace 89 \rbrace & &\lbrace 94\rbrace\to\lbrace 97 \rbrace & &\lbrace 94\rbrace\to\lbrace 99 \rbrace\\
    &78: 1,2,5,6 & &78: \text{n/a} & &78: 4 & &78: 3          & &78: \text{n/a} & &78: \text{n/a}\\
    &82: 1,2,5,6 & &82: \text{n/a} & &82: 4 & &82: 3          & &82: \text{n/a} & &82: \text{n/a}\\
    &86: 1,2,5,6 & &86: \text{n/a} & &86: 4 & &86: 3          & &86: \text{n/a} & &86: \text{n/a}\\
    &90: 1,2,5,6 & &90: \text{n/a} & &90: 4 & &90: 3          & &90: \text{n/a} & &90: \text{n/a}\\
    &94: 1,2,6   & &94: \text{n/a} & &94: 4 & &94: \text{n/a} & &94: 3          & &94: 5\\
\end{align*}
Notice, four of the five squares in $\lbrace 94 \rbrace$ are such that when the player rolls a non-4, four of the possible five rolls will send the player to $\lbrace w \rbrace$ (column 1 above, 78-90), and the remaining square in $\lbrace 94 \rbrace$ is such that when the player rolls a non-4, three of the possible five rolls will send the player to $\lbrace w \rbrace$ (column 1 above, 94). Therefore,
\begin{align*}
    P\left(\lbrace 94\rbrace \to \lbrace w \rbrace\right) &= \frac{4}{5}\cdot\frac{4}{5}[1-P(4)] + \frac{1}{5}\cdot\frac{3}{5}[1-P(4)]\\
    &=\frac{16}{25}\cdot[1-P(4)] + \frac{3}{25}\cdot[1-P(4)]\\
    &= \frac{19}{25}\cdot[1-P(4)]
\end{align*}
Similarly,
\begin{align*}
    P(\lbrace 94 \rbrace \to \lbrace 60 \rbrace) &= 0\\\\
    P(\lbrace 94 \rbrace \to \lbrace 94 \rbrace) &= \frac{5}{5}\cdot P(4)\\ &= P(4)\\\\
    P(\lbrace 94 \rbrace \to \lbrace 89 \rbrace) &= \frac{4}{5}\cdot\frac{1}{5}[1-P(4)]\\
    &= \frac{4}{25}[1-P(4)]\\\\
    P(\lbrace 94 \rbrace \to \lbrace 97 \rbrace) &= \frac{1}{5}\cdot \frac{1}{5} \cdot [1-P(4)]\\
    &= \frac{1}{25}[1-P(4)]\\\\
    P(\lbrace 94 \rbrace \to \lbrace 99 \rbrace) &=\frac{1}{5}\cdot \frac{1}{5} \cdot [1-P(4)]\\
    &= \frac{1}{25}[1-P(4)]
\end{align*}
Repeating this process for all five other possible starting states $\left( \lbrace w \rbrace, \lbrace 60 \rbrace, \lbrace 89 \rbrace, \lbrace 97 \rbrace, \text{ and } \lbrace 99 \rbrace \right)$ we find the transition probabilities necessary to construct the Markov model in figure \ref{tp_4} for the player's transition between the outlined 4-win path and all 4-loops in figure \ref{4class}. We omit writing all the transition possibilities from each square in the other 4-loops; it's the same process shown for $\lbrace 94 \rbrace$, and trust us, figure \ref{tp_4} is much more interesting to look at.

\newpage

Next, with figure \ref{tp_4} in hand, we find the proportion of turns the player spends in all possible states by solving the following system of equations, where $P(\lbrace i\rbrace \to \lbrace j \rbrace)$ are given in figure \ref{tp_4}.
\begin{align*}
    p_w &= p_w\cdot P(\lbrace w \rbrace \to \lbrace w \rbrace) + p_{60}\cdot P(\lbrace 60 \rbrace \to \lbrace w \rbrace) + p_{89}\cdot P(\lbrace 89 \rbrace \to \lbrace w \rbrace)\\ &+ p_{94}\cdot P(\lbrace 94 \rbrace \to \lbrace w \rbrace) + p_{97}\cdot P(\lbrace 97 \rbrace \to \lbrace w \rbrace) + p_{99}\cdot P(\lbrace 99 \rbrace \to \lbrace w \rbrace)\\\\
    p_{60} &= p_w\cdot P(\lbrace w \rbrace \to \lbrace 60 \rbrace) + p_{60}\cdot P(\lbrace 60 \rbrace \to \lbrace 60 \rbrace) + p_{89}\cdot P(\lbrace 89 \rbrace \to \lbrace 60 \rbrace)\\ &+ p_{94}\cdot P(\lbrace 94 \rbrace \to \lbrace 60 \rbrace) + p_{97}\cdot P(\lbrace 97 \rbrace \to \lbrace 60 \rbrace) + p_{99}\cdot P(\lbrace 99 \rbrace \to \lbrace 60 \rbrace)\\\\
    p_{89} &= p_w\cdot P(\lbrace w \rbrace \to \lbrace 89 \rbrace) + p_{60}\cdot P(\lbrace 60 \rbrace \to \lbrace 89 \rbrace) + p_{89}\cdot P(\lbrace 89 \rbrace \to \lbrace 89 \rbrace)\\ &+ p_{94}\cdot P(\lbrace 94 \rbrace \to \lbrace 89 \rbrace) + p_{97}\cdot P(\lbrace 97 \rbrace \to \lbrace 89 \rbrace) + p_{99}\cdot P(\lbrace 99 \rbrace \to \lbrace 89 \rbrace)\\\\
    p_{94} &= p_w\cdot P(\lbrace w \rbrace \to \lbrace 94 \rbrace) + p_{60}\cdot P(\lbrace 60 \rbrace \to \lbrace 94 \rbrace) + p_{89}\cdot P(\lbrace 89 \rbrace \to \lbrace 94 \rbrace)\\ &+ p_{94}\cdot P(\lbrace 94 \rbrace \to \lbrace 94 \rbrace) + p_{97}\cdot P(\lbrace 97 \rbrace \to \lbrace 94 \rbrace) + p_{99}\cdot P(\lbrace 99 \rbrace \to \lbrace 94 \rbrace)\\\\
    p_{97} &= p_w\cdot P(\lbrace w \rbrace \to \lbrace 97 \rbrace) + p_{60}\cdot P(\lbrace 60 \rbrace \to \lbrace 97 \rbrace) + p_{89}\cdot P(\lbrace 89 \rbrace \to \lbrace 97 \rbrace)\\ &+ p_{94}\cdot P(\lbrace 94 \rbrace \to \lbrace 97 \rbrace) + p_{97}\cdot P(\lbrace 97 \rbrace \to \lbrace 97 \rbrace) + p_{99}\cdot P(\lbrace 99 \rbrace \to \lbrace 97 \rbrace)\\\\
    p_{99} &= p_w\cdot P(\lbrace w \rbrace \to \lbrace 99 \rbrace) + p_{60}\cdot P(\lbrace 60 \rbrace \to \lbrace 99 \rbrace) + p_{89}\cdot P(\lbrace 89 \rbrace \to \lbrace 99 \rbrace)\\ &+ p_{94}\cdot P(\lbrace 94 \rbrace \to \lbrace 99 \rbrace) + p_{97}\cdot P(\lbrace 97 \rbrace \to \lbrace 99 \rbrace) + p_{99}\cdot P(\lbrace 99 \rbrace \to \lbrace 99 \rbrace)\\\\
    1 &= p_w + p_{60} + p_{89} + p_{94} + p_{97} + p_{99}
\end{align*}
Feeding this to CoCal SageMath, we calculate,

\begin{figure}[H]
    \centering
    \caption{Proportion of Turns Occupying Each State in Fig. \ref{tp_4}}
    \label{prop1}
    \vspace*{0.5cm}
    \includegraphics[scale=1]{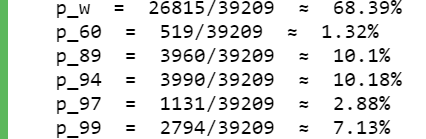}
\end{figure}
\noindent Similar to section \ref{on5}, the proportion of turns spent in each state is independent of $P(4)$.

\begin{figure}[H]
    \caption{Markov Chain for the Player’s Transition Between the 4-Win Path and 4-loops}
    \label{tp_4s}
    \begin{center}
        \begin{tikzpicture}
                % Add the states
                \node[state] at (-3,0)       (w) {Win};
                \node[state] at (3,0)       (l) {Loop};

                \draw[every loop]
                    (w) edge[bend left, auto=left] node {$\frac{37}{155} \cdot [1-P(4)]$}(l)
                    (w) edge[loop left] node {$P(4) + \frac{118}{155} \cdot [1-P(4)]$}(w)
                    (l) edge[bend left, auto=left] node {$\frac{6401}{12394} \cdot [1-P(4)]$}(w)
                    (l) edge[loop right] node {$P(4) + \frac{5993}{12394} \cdot [1-P(4)]$}(l);
             \end{tikzpicture}
    \end{center}
\end{figure}
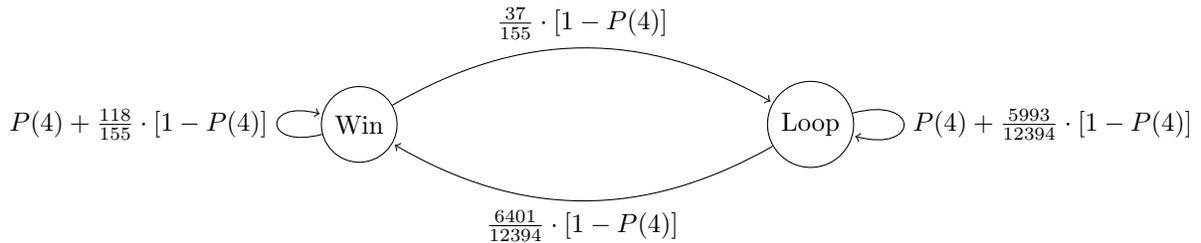

Yet, in order to derive a function of $x$ for the expected duration of each game, we need to know the probability the player transitions from any of the loops to the win state. To find this, we will simplify the model down to only two states, just as we had in figure \ref{tp_5}. Let $\Omega$ be the set of all loop states. The probability the player transitions from the win state to any loop is equal to the sum of all $P(Win \to j)$ for $j \in \Omega$. Therefore,
\[P\left(\lbrace w \rbrace \to \Omega\right) = \frac{37}{155}[1-P(4)]\]
Furthermore, $P\left(\lbrace w \rbrace \to \lbrace w \rbrace\right)$ remains as shown in figure \ref{tp_4}. Next, to find $P\left( \Omega \to \lbrace w \rbrace \right)$, we are unfortunately unable sum all $P(i \to Win)$ for $i \in \Omega$ since this will clearly be greater than 1. So, we need to scale these transition probabilities down. From figure \ref{prop1}, the probability the player is in any loop is
\[p_{60} + p_{89} + p_{94} + p_{97} + p_{99} = \frac{12394}{39209}\]
So, we need $X_i$ the probability the player is in a specific loop state $i$, given they are in any loop state. That is, with $p_i$ the probability the player is in the loop $i \in \Omega$,
\begin{align*}
    p_i &= \frac{12394}{39209} \cdot X_i\\
    \implies X_i &= p_i \cdot \frac{39209}{12394}\\
    &= \frac{n(p_i)}{12394}
\end{align*}
where $n(p_i)$ is the numerator of $p_i$, since the denominator of $p_i$ is $39209$ for all $i$ (figure \ref{prop1}). So, $n(p_i) / 12394$ is the coefficient of $P(i \to Win)$ that ensures \[ \frac{n(p_i)}{12394} \sum_{i \in \Omega} P\left( i \to \lbrace w \rbrace\right) = 1\]
So,
\begin{align*}
    P(Loop \to Win) &= \sum_{i\in\Omega} P\left(i \to \lbrace w \rbrace\right) \cdot X_i\\
    &= [1-P(4)]\cdot\left(\frac{2}{5}\cdot\frac{n(p_{60})}{12394} + \frac{3}{5}\cdot\frac{n(p_{89})}{12394} + \frac{19}{25}\cdot\frac{n(p_{94})}{12394} + \frac{1}{5}\cdot\frac{n(p_{97})}{12394} +\frac{1}{5}\cdot\frac{n(p_{99})}{12394}\right)\\
    &= [1-P(4)]\cdot\left( \frac{2}{5}\cdot\frac{519}{12394} + \frac{3}{5}\cdot\frac{3960}{12394} + \frac{19}{25}\cdot\frac{3990}{12394} + \frac{1}{5}\cdot\frac{1131}{12394} + \frac{1}{5}\cdot\frac{2794}{12394} \right)\\
    &= [1-P(4)]\cdot\left(\frac{6401}{12394}\right)\\
    P(Loop \to Loop) &= \sum_{i\in\Omega} \left(\sum_{j \in \Omega} P(i \to j) \cdot X_i\right)\\
    &= P(4)\cdot\left( \sum_{i\in\Omega}\frac{n(p_i)}{12394} \right) + [1-P(4)]\cdot \left( \left[ \frac{1}{5} + \frac{2}{5} \right] \cdot \frac{n(p_{60})}{12394} + \left[ \frac{2}{5} \right] \cdot \frac{n(p_{89})}{12394}\right.\\ 
    &+ \left.\left[ \frac{4}{25} + \frac{1}{25} + \frac{1}{25} \right] \cdot \frac{n(p_{94})}{12394} + \left[ \frac{1}{5} + \frac{1}{5} + \frac{2}{5} \right] \cdot \frac{n(p_{97})}{12394}  + \left[ \frac{4}{5} \right]\cdot \frac{n(p_{99})}{12394}\right)\\
    &= P(4) \cdot 1 + [1-P(4)]\cdot \left( \left[\frac{3}{5}\right]\cdot \frac{519}{12394} + \left[ \frac{2}{5}\right]\cdot \frac{3960}{12394} + \left[ \frac{6}{25} \right]\cdot \frac{3990}{12394} + \left[ \frac{4}{5} \right]\cdot\frac{1131}{12394} + \left[\frac{4}{5}\right]\cdot\frac{2794}{12394} \right)\\
    &= P(4) + [1-P(4)]\cdot\frac{5993}{12394}
\end{align*}
Finally, we construct the Markov model in figure \ref{tp_4s}.
We are now able to construct a function of $x$, where $x$ is the average number of turns per non-4 roll, that approximates the expected duration of the game. For $P(4) = p(x) = \frac{x-1}{x}$, the probability the game ends in exactly 31 turns is $p(x)^{31}$ and this contributes $31$ moves to the expected value. The probability this does not happen is $1-p(x)^{31}$. In this case, the player transitions between the 4-win path and 4-loops roughly according to the Markov model in figure \ref{tp_4s}. According to this model, the probability the player transitions from any loop state to the 4-win path in one turn is $\frac{6401}{12394}\cdot [1-P(4)]$. Therefore, the expected number of times the player needs to roll a non-4 to escape is $\frac{12394}{6401}$. Yet, they only need to roll a non-4 for as long as they are in a loop state and the proportion of turns the player spends in loops states is $\frac{12394}{39209}$. Therefore,
\begin{align*}
    E_4(x) &= p(x)^{31} \cdot 31 + \left[1-p(x)^{31}\right]\cdot\frac{12394}{39209}\cdot\frac{12394}{6401}x\\
    &\approx p(x)^{31} \cdot 31 + \left[1-p(x)^{31}\right]\cdot 0.612x
\end{align*}
roughly approximates the expected number of turns per game as a function of the average number of rolls per single non-4 roll. Notice,
\[\sigma \approx \lim_{x \to \infty}E_4(x) \approx 49.97\]
Clearly, $49.94 \not = 46.98$; that is, our approximation of $\sigma$ from $E_4$ does not equal true $\sigma$. Notice, for the case of $P(5) \to 1$ in section \ref{on5}, the board is fairly simple and the simplified model in figure \ref{tp_5} is within 0.16\% of true $\sigma$. However, for the case of $P(4)\to 1$, the board is much more complicated; the on ramp to the outlined 4-win path in figure \ref{4class} is much larger than the on ramp to the outlined 5-win path in figure \ref{5class}. Since neither simplified model considers cases where the player rolls a non-$\delta$ while occupying any square in any on ramp, simplifying the model has much greater effect on its accuracy for the case of $P(4) \to 1$ than it did for the case of $P(5) \to 1$. This is likely a significant contributing factor to why the model in figure \ref{tp_4s} is only within 6.5\% of true $\sigma$ whereas the model in figure \ref{tp_5} is within 0.16\%.

\subsection{Open Question}\label{openQ1}
In both section \ref{on5} and section \ref{on4} we derived a function of $x$, the average number of rolls per single non-$\delta$ roll, for $\delta \in \lbrace 4,5\rbrace$, whose limit as $x\to \infty$ was the expected duration of the game as $P(\delta) \to 1$. Let $\psi$ be the number of turns per game where $P(\delta) = 1$. Then, since $\sigma$ depends on $x$, we found previously
\[\sigma(x) \approx E_\delta(x) = p(x)^\psi\cdot \psi +\left[1-p(x)^\psi\right]\cdot cx\]
Where $P(\delta) = p(x) = \frac{x-1}{x}$. I claim, if $\sigma(x)$ is known, we may find a better coefficient $c$ of $x$ such that $\lim_{x \to \infty} E_\delta(x)$ will be within any tolerance of $\sigma(x)$; dependent on only how many decimals we care to write. Notice,
\begin{align*}
\sigma(x) = \frac{\text{Number of Turns}}{\text{Number of Games}} &= \frac{(x-1) + cx}{(x-1)/\psi + 1}
\end{align*}
That is, with $P(\delta)$ close to $1$, on average, of the first $x$ rolls, $x-1$ of them are rolls of $\delta$; which yield $\frac{x-1}{\psi}$ games since $\psi$ is the number turns per game if all turns are rolls of $\delta$. Then, the game where the player rolls their first non-$\delta$ lasts $cx$ turns. Thus, the total number of turns is $x-1 + cx$ and the total number of games is $\frac{x-1}{\psi} + 1$. And so,
\begin{align*}
   \sigma(x) &= \frac{(x-1) + cx}{(x-1)/\psi + 1}\\
   &= (x-1+cx)\cdot\frac{1}{(x-1)/\psi + 1}\\
   &= (x-1+cx)\cdot\frac{\psi}{x-1 + \psi}\\
   &= \frac{(\psi+c\psi)x - \psi}{x-1 + \psi}
\end{align*}
Yet, the above equation requires $P(\delta) \to 1$ which happens as $x \to \infty$. So,
\begin{align*}
    \lim_{x \to \infty} \sigma(x) &= \lim_{x \to \infty} \frac{(\psi+c\psi)x - \psi}{x-1 + \psi}\\
    \sigma &= \psi+c\psi\\
    &= (1+c)\psi\\
    \implies c &= \frac{\sigma}{\psi}-1
\end{align*}
Therefore,
\begin{align*}
    c_{\delta = 5} &= \frac{82.37}{16} - 1 & c_{\delta = 4} &= \frac{46.98}{31} - 1\\
    &= \frac{66.37}{16} & &=  \frac{15.98}{31}\\
    &\approx 4.148   & &\approx 0.515
\end{align*}
And using $c_\delta$ in $E_\delta$,
\begin{align*}
    \lim_{x \to \infty} E_5(x) &= 82.37 = \sigma & \lim_{x \to \infty} E_4(x) &= 46.98 = \sigma
\end{align*}
Whereas, when we used the simplified Markov model of the game we approximated,
\begin{align*}
    c_5 &= \frac{77}{93} \cdot 5 & c_4 &=\frac{12394}{39209}\cdot\frac{12394}{6401}\\
    &\approx 4.140 & &\approx 0.612
\end{align*}
and
\begin{align*}
    \lim_{x \to \infty} E_5(x) &\approx 82.24 & \lim_{x \to \infty} E_4(x) &\approx 49.97
\end{align*}
The main difference being, $c = \sigma / \psi - 1$ gives us an exact value for $c$ but requires the value of $\sigma$ be known; whereas, the simplified Markov model of the game lets us approximate $c$ when $\sigma$ is unknown. This begs the question, is it possible to approximate the value of $c$ better than we did in section \ref{on5} and \ref{on4} such that $\sigma$ is not involved in the approximation?

\section{Introducing Strategy}\label{IntroStrat}
Suppose a fair coin is added to the base game as a method of changing the player's position on the board in addition to the fair die in the following fashion. Each turn, the player may choose whether or not they flip the coin after they roll the fair die and advance their position but before they ascend or descend any ladder or chute. Should they choose not to flip the coin, the player ascends or descends a ladder or chute if they landed on one and their turn is over. However, if they choose to flip the coin and they get heads then they advance one square, and if they get tails then they move back one square; after the coin flip and resulting position change, they will ascend or descend any ladder or chute if they landed on one and their turn is over. For example, suppose the player is on square-53 and they roll a three. They advance to square-56, which is at the top of a slide that leads back to square-53. They do not immediately descend the slide. The player elects to flip the coin and they get tails. They are now on square-55, which is not at the bottom of a ladder or top of a slide, so their turn is over. Clearly, to never flip the coin is to play the base game as is; label this Strategy 0. 
\subsection{Strategy Descriptions}\label{Strat Desc}
The strategy the player employs to determine whether or not to flip the coin at each turn may have a significant impact on the average duration of the game. Below we define seven different strategies that we will later analyze. Ladders are shown in green while chutes are red. A square colored \includegraphics[scale=0.25]{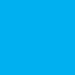} indicates the player will choose to flip the coin when they land on that square, whereas a square colored \includegraphics[scale=0.25]{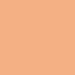} indicates the player will choose not to flip the coin when they land on that square, all according the specified strategy.
\begin{itemize}
    \item Strategy -1: Randomly flip the coin.
    \item Strategy 0: Never flip the coin; this is the base game.
    
    \begin{figure}[H]
        \centering
        \caption{Strategy 0 Square Classification}
        \includegraphics[scale=0.5]{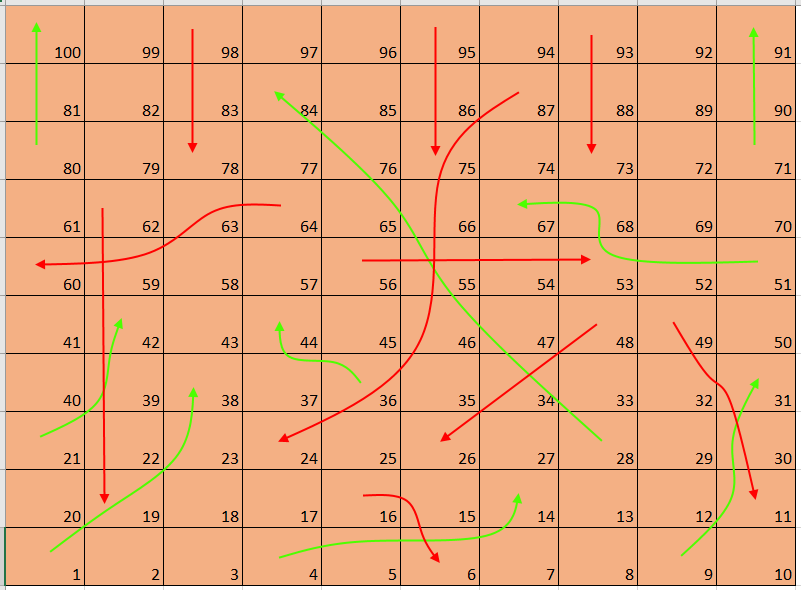}
        \label{strat0flips}
    \end{figure}

    \item Strategy 1: Do not flip the coin when on square-100; otherwise randomly flip the coin.
    \item Strategy 2: Always flip the coin, unless on square-100.
    
    \begin{figure}[H]
        \centering
        \caption{Strategy 2 Square Classification}
        \includegraphics[scale=0.5]{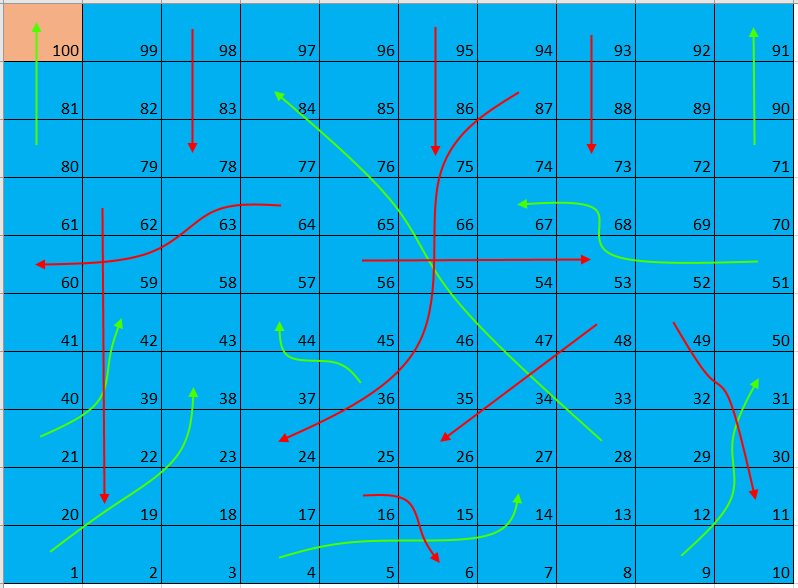}
        \label{strat2flips}
    \end{figure}
    \vspace*{2cm}

\item Strategy 3: Always flip the coin, unless at the base of a ladder or on square-100.
\begin{figure}[H]
    \centering
    \caption{Strategy 3 Square Classification}
    \includegraphics[scale=0.5]{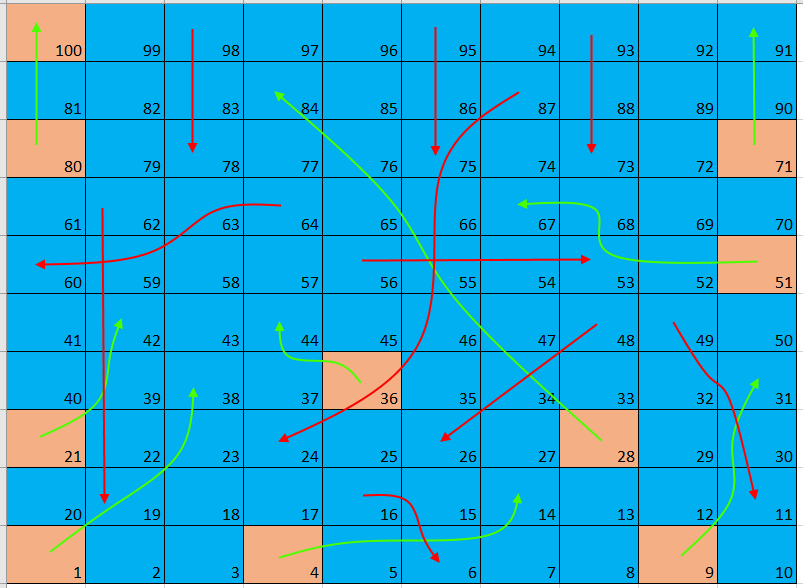}
    \label{strat3flips}
\end{figure}

\item Strategy 4: Only flip the coin when at the top of a chute.
\begin{figure}[H]
    \centering
    \caption{Strategy 4 Square Classification}
    \includegraphics[scale=0.5]{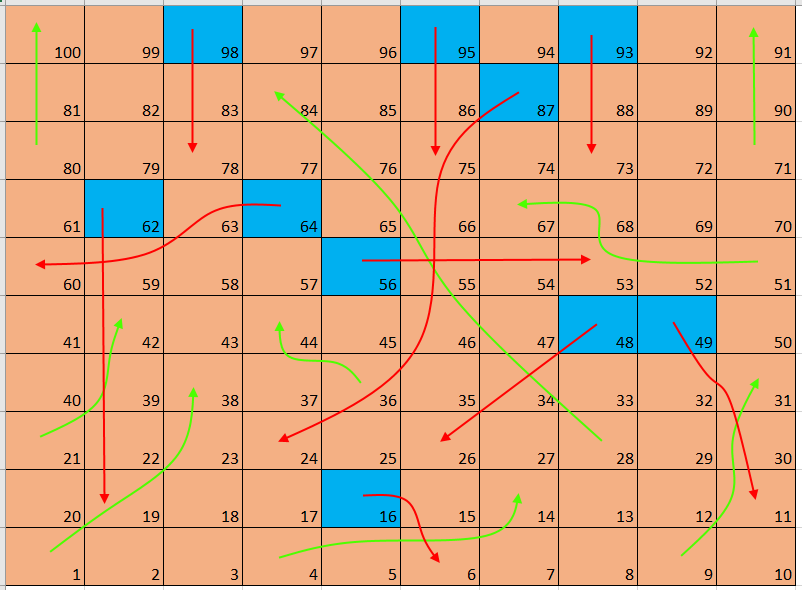}
    \label{strat4flips}
\end{figure}

\vspace*{2cm}

\item Strategy 5: Only flip the coin when at the top of a chute, or when one square away from the base of any ladder.
\begin{figure}[H]
    \centering
    \caption{Strategy 5 Square Classification}
    \includegraphics[scale=0.5]{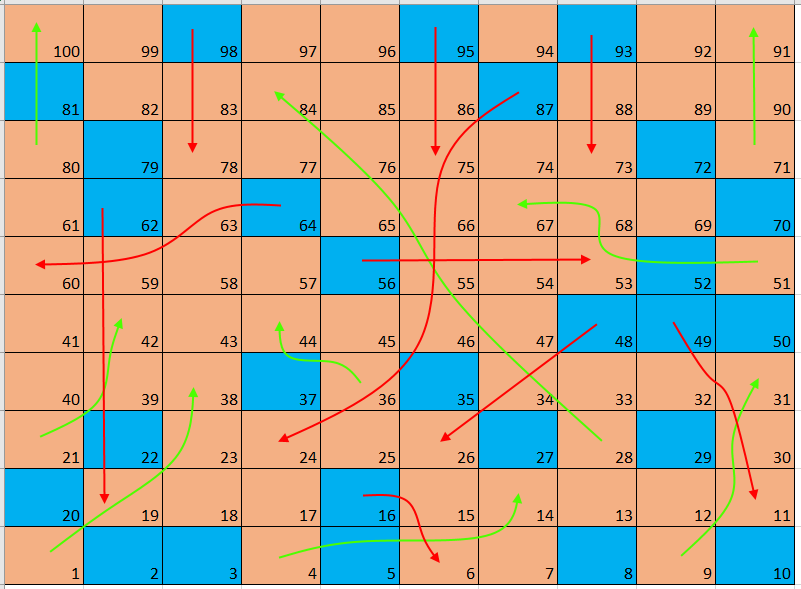}
    \label{strat5flips}
\end{figure}

\item Strategy 6: Never flip the coin when at the base of a ladder, or when one square from the top of a chute, or when on square-100. Otherwise flip.
\begin{figure}[H]
    \centering
    \caption{Strategy 6 Square Classification}
    \includegraphics[scale=0.5]{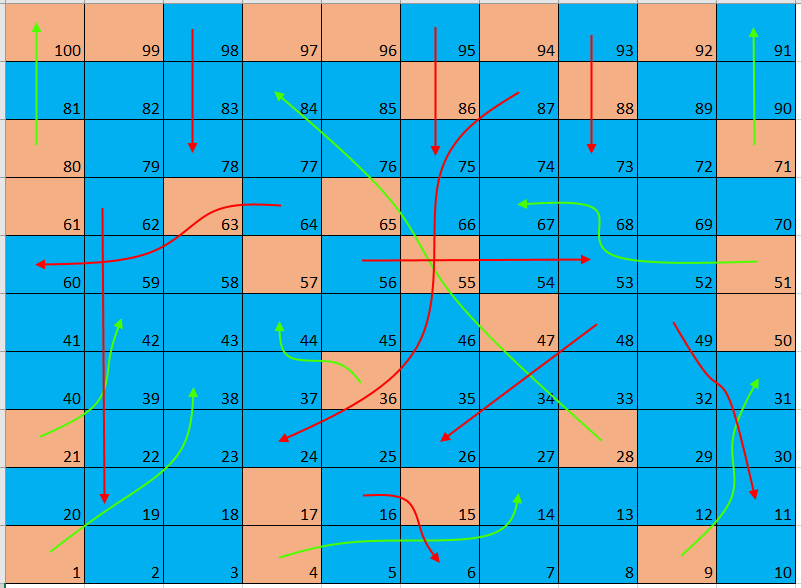}
    \label{strat6flips}
\end{figure}
\end{itemize}

\subsection{Determining $\sigma$ of Each Strategy}\label{sratsigma}
We utilized Python to simulate many games to approximate the average number of turns per game for each strategy. The following table shows the calculated $\sigma$ for each strategy after the specified number of simulated games.
\begin{figure}[H]
    \centering
    \caption{Average Duration Over Simulated Games Per Strategy}
    \vspace*{0.25cm}
    \includegraphics[scale=1.25]{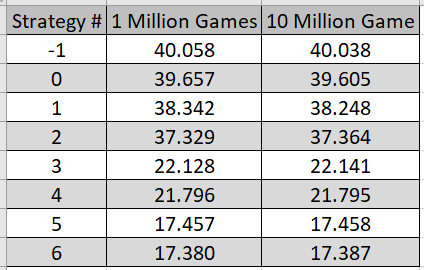}
    \label{fig:my_label}
\end{figure}

\noindent Interestingly, Strategy 4 is very similar to, but not exactly like, playing the base game with only a fair die and on the same board but with all the chutes removed. After simulating ten million games of this description, we calculated $\sigma$ to be $21.843$. A significant difference between a game of this description and Strategy 4 is that square-48 and square-49 are both at the top of a chute; so, when the player lands on either square and flips the coin it is still possible they will descend a chute. Whereas, if all chutes are removed from the board then it is clearly impossible to descend any chute. This seems to suggest that Strategy 4 should yield a larger $\sigma$ than the game described, yet it does not; Strategy 4 appears to yield a lower $\sigma$.

\newpage

\section{SageMath \& Python Code}\label{Code}
\noindent Python Code: Simulate Games Optional Weighted Die............................................\ref{py1}\\
\noindent SageMath Code: Markov Model of Chutes \& Ladders..............................................\ref{sm1}\\
\noindent Python Code: Simulate All Start Positions Over All $P(\delta)=1$..................................\ref{py2}\\
\noindent Python Code: Simulate Game With Coin Strategies.................................................\ref{py3}

\begin{center}
\lstset{
    language=python,
    tabsize=2,
    %frame=lines,
    frame=shadowbox,
    rulesepcolor=\color{gray},
    xleftmargin=-50pt,
    xrightmargin=-50pt,
    framexleftmargin=15pt,
    keywordstyle=\color{blue},
    commentstyle=\color{60color},
    stringstyle=\color{red},
    numbers=left,
    numberstyle=\tiny,
    numbersep=5pt,
    breaklines=true,
    showstringspaces=false,
    basicstyle=\footnotesize,
    emph={str},emphstyle={\color{magenta}},
    caption={Python Code: Simulate Games Optional Weighted Die}, captionpos=t
    }
    
\label{py1}
\begin{lstlisting}
import random


def Welcome():
    print("Welcome to the Chutes & Ladders simulation!")
    Start()

def Start():
    #initialize vars
    nums = [1,2,3,4,5,6]
    tbt_bool = True
    sum_bool = False
    weight_bool = False
    total_games = ""
    weight_num = 0
    weight_amt = 0
    w=[0]
    p=0

    #input
    tbt_YN = input("Do you want a turn by turn readout? (Only able to simulate a game at a time) ")
    if tbt_YN == "Y" or tbt_YN == "y":
        tbt_bool = True
    elif tbt_YN == "N" or tbt_YN == "n":
        tbt_bool = False
    else:
        print("Only 'y' or 'n'.")
        Start()
    if tbt_bool == False:
        sum_YN = input("Do you want a summary readout for each game? ")
        if sum_YN == "Y" or sum_YN == "y":
            sum_bool = True
        elif sum_YN == "N" or sum_YN == "n":
            sum_bool = False
        else:
            print("Only 'y' or 'n'.")
            Start()

    if tbt_bool == True:
        total_games = 1
    else:
        total_games = input("Number of games to simulate: ")
        try:
            int(total_games)
        except:
            print("You must enter a positive integer.")
            Start()
        if int(total_games) <= 0:
            print("You must enter a positive integer.")
            Start()
    weight_YN = input("Do you want to simulate a weighted die? ")
    if weight_YN == "Y" or weight_YN == "y":
        weight_bool = True
        weight_num = input("Which number would you like to weight? (1-6) ")
        try:
            int(weight_num)
        except:
            print("Must enter an integer between 1 and 6. A")
            Start()
    elif weight_YN == "N" or weight_YN == "n":
        weight_bool = False
    else:
        print("Only 'y' or 'n'.")
        Start()
    if nums.count(int(weight_num)) == 0 and weight_bool == True:
        print("Must enter an integer between 1 and 6. B")
        Start()
    if nums.count(int(weight_num)) > 0 and weight_bool == True:
        print("Set probability of",weight_num,"appearing: ", end='')
        weight_amt = input("")
        try:
            float(weight_amt)
        except:
            print("Must be number between 0 and 1")
            Start()
        if float(weight_amt) < 0 or float(weight_amt) > 1:
            print("Must be number between 0 and 1")
            Start()


    if weight_bool == True:
        nums.remove(int(weight_num))
        weight_cmp = round(1 - float(weight_amt),6)
        for i in range(1,7,1):
            if nums.count(i) > 0:
                p = p + round(weight_cmp * 2000000,0)
                w.append(p)
            else:
                p = p + round(float(weight_amt) * 10000000,0)
                w.append(p)
        print(w)


    Game(tbt_bool,sum_bool,int(total_games), w)



def dist(w):
    broll = random.randint(1,10000000)
    if broll <= w[1]:
        return int(1)
    elif broll <= w[2]:
        return int(2)
    elif broll <= w[3]:
        return int(3)
    elif broll <= w[4]:
        return int(4)
    elif broll <= w[5]:
        return int(5)
    elif broll <= w[6]:
        return int(6)





def Game(display, summary, total_games, w):
    game_num = 1
    notified = False
    long_counter = 0
    # Set up the boar#
    #Pos:       1   2  3  4   5  6  7  8  9  10 11 12 13 14 15 16 17 18 19 20  21 22 23 24 25 26 27  28 29 30 31 32 33 34 35  36 37 38 39 40 41 42 43 44 45 46 47  48  49 50
    board = [0, 38, 0, 0, 14, 0, 0, 0, 0, 31, 0, 0, 0, 0, 0, 0, 6, 0, 0, 0, 0, 42, 0, 0, 0, 0, 0, 0, 84, 0, 0, 0, 0, 0, 0, 0, 44, 0, 0, 0, 0, 0, 0, 0, 0, 0, 0, 0, 26, 11, 0,
                67, 0, 0, 0, 0, 53, 0, 0, 0, 0, 0, 19, 0, 60, 0, 0, 0, 0, 0, 0, 91, 0, 0, 0, 0, 0, 0, 0, 0, 100, 0, 0, 0, 0, 0, 0, 24, 0, 0, 0, 0, 0, 73, 0, 75, 0, 0, 78, 0, 0]
    #Pos:       51 52 53 54 55  56 57 58 59 60 61  62 63  64 65 66 67 68 69 70  71 72 73 74 75 76 77 78 79   80 81 82 83 84 85 86  87 88 89 90 91 92  93 94  95 96 97  98 99 100

    # Initialize player position
    player_pos = 0

    #Initialize move counter
    move_number = 0
    #Initialize game status
    gameOn = True

    #Initialize total move record
    moves_per_game = []
    #Initialize total moves
    total_moves = 0




    while gameOn == True:
        move_number += 1
        #1st Game
        if game_num == 1 and move_number == 1 and summary == True:
            print("Game #1!")
        # Roll the dice
        if display == True:
            print("Move #",move_number,":", sep='')




        #Roll
        if len(w) == 1:
            roll = random.randint(1, 6)
        else:
            roll = dist(w)
        if display == True:
            print("You rolled a", roll)




        #Long games
        if move_number > 1000 and notified == False:
            notified = True
            #print("Game", game_num, "passed 1k moves")
            #break


        # Move the player
        player_pos += roll


        #Check new location
        if player_pos > 100:
            player_pos -= roll
            if display == True:
                print("You must roll exactly to land on final square. Roll again.")
                print("You have returned to square", player_pos)
        else:
            if board[player_pos] != 0:
                if display == True:
                    print("You landed on a", "chute" if board[player_pos] < player_pos else "ladder!")
                player_pos = board[player_pos]
            if display == True:
                print("You are now on square", player_pos)
        if player_pos == 100:
            # WIINER WINNER CHICEKN DINNER
            moves_per_game.append(move_number)
            if display == True:
                print("You won!")
            if summary == True:
                print("Total number of moves:", move_number)
            if game_num == total_games:
                gameOn = False
                for i in range(total_games):
                    total_moves += moves_per_game[i-1]
                print("___________________________________________")
                print("Simulation comeplete. Played ", game_num, " game(s). Total moves: ", total_moves, ". Average moves per game: ", total_moves/total_games,".", sep='')
                print("Number of games over 1k moves", long_counter)
            else:
                if notified == True:
                    #print("Game", game_num, "took", move_number, "moves.")
                    long_counter += 1
                game_num += 1
                move_number = 0
                player_pos = 0
                notified = False
                if summary == True:
                    print("___________________________________________")
                    print("NEW GAME! This is now game #",game_num,".",sep='')









Welcome()
\end{lstlisting}
\end{center}

\newpage

\begin{center}
\lstset{
    language=python,
    tabsize=1,
    frame=none,
    %frame=shadowbox,
    rulesepcolor=\color{gray},
    xleftmargin=-25pt,
    %xrightmargin=-50pt,
    %framexleftmargin=15pt,
    keywordstyle=\color{blue},
    commentstyle=\color{60color},
    basicstyle=\tiny,
    stringstyle=\color{red},
    numbers=none,
    %numberstyle=\tiny,
    %numbersep=5pt,
    breaklines=true,
    %postbreak=\mbox{\textcolor{red}{$\hookrightarrow$}\space},
    showstringspaces=false,
    emph={str},emphstyle={\color{magenta}},
    captionpos=t, caption={SageMath Code: Markov Model of Chutes \& Ladders}
    }
\label{sm1}
\begin{lstlisting}
var('o, t, r, f, v, s')
o = 1/6
t = 1/6
r = 1/6
f = 1/6
v = 1/6
s = 1/6
spos=vector(
[1,0,0,0,0,0,0,0,0,0,0,0,0,0,0,0,0,0,0,0,0,0,0,0,0,0,0,0,0,0,0,0,0,0,0,0,0,0,0,0,0,
0,0,0,0,0,0,0,0,0,0,0,0,0,0,0,0,0,0,0,0,0,0,0,0,0,0,0,0,0,0,0,0,0,0,0,0,0,0,0,0]
)

one=vector(
[1,1,1,1,1,1,1,1,1,1,1,1,1,1,1,1,1,1,1,1,1,1,1,1,1,1,1,1,1,1,1,1,1,1,1,1,1,1,1,1,1,
1,1,1,1,1,1,1,1,1,1,1,1,1,1,1,1,1,1,1,1,1,1,1,1,1,1,1,1,1,1,1,1,1,1,1,1,1,1,1,1]
)

###########################################################################################################################

T=matrix([
[0,t,r,v,s,0,0,0,0,0,0,f,0,0,0,0,0,0,0,0,0,0,0,0,0,0,0,0,0,0,0,o,0,0,0,0,0,0,0,0,0,
0,0,0,0,0,0,0,0,0,0,0,0,0,0,0,0,0,0,0,0,0,0,0,0,0,0,0,0,0,0,0,0,0,0,0,0,0,0,0,0],#0
[0,0,o,r,f,v,s,0,0,0,0,t,0,0,0,0,0,0,0,0,0,0,0,0,0,0,0,0,0,0,0,0,0,0,0,0,0,0,0,0,0,
0,0,0,0,0,0,0,0,0,0,0,0,0,0,0,0,0,0,0,0,0,0,0,0,0,0,0,0,0,0,0,0,0,0,0,0,0,0,0,0],#2
[0,0,0,t,r,f,v,0,0,0,0,o,0,0,0,0,0,0,0,0,0,0,0,0,0,s,0,0,0,0,0,0,0,0,0,0,0,0,0,0,0,
0,0,0,0,0,0,0,0,0,0,0,0,0,0,0,0,0,0,0,0,0,0,0,0,0,0,0,0,0,0,0,0,0,0,0,0,0,0,0,0],#3
[0,0,0,0,o,t,r,v,s,0,0,0,0,0,0,0,0,0,0,0,0,0,0,0,0,f,0,0,0,0,0,0,0,0,0,0,0,0,0,0,0,
0,0,0,0,0,0,0,0,0,0,0,0,0,0,0,0,0,0,0,0,0,0,0,0,0,0,0,0,0,0,0,0,0,0,0,0,0,0,0,0],#5
[0,0,0,0,0,o,t,f,v,s,0,0,0,0,0,0,0,0,0,0,0,0,0,0,0,r,0,0,0,0,0,0,0,0,0,0,0,0,0,0,0,
0,0,0,0,0,0,0,0,0,0,0,0,0,0,0,0,0,0,0,0,0,0,0,0,0,0,0,0,0,0,0,0,0,0,0,0,0,0,0,0],#6
[0,0,0,0,0,0,o,r,f,v,s,0,0,0,0,0,0,0,0,0,0,0,0,0,0,t,0,0,0,0,0,0,0,0,0,0,0,0,0,0,0,
0,0,0,0,0,0,0,0,0,0,0,0,0,0,0,0,0,0,0,0,0,0,0,0,0,0,0,0,0,0,0,0,0,0,0,0,0,0,0,0],#7
[0,0,0,0,0,0,0,t,r,f,v,s,0,0,0,0,0,0,0,0,0,0,0,0,0,o,0,0,0,0,0,0,0,0,0,0,0,0,0,0,0,
0,0,0,0,0,0,0,0,0,0,0,0,0,0,0,0,0,0,0,0,0,0,0,0,0,0,0,0,0,0,0,0,0,0,0,0,0,0,0,0],#8
[0,0,0,0,s,0,0,0,o,t,r,f,v,0,0,0,0,0,0,0,0,0,0,0,0,0,0,0,0,0,0,0,0,0,0,0,0,0,0,0,0,
0,0,0,0,0,0,0,0,0,0,0,0,0,0,0,0,0,0,0,0,0,0,0,0,0,0,0,0,0,0,0,0,0,0,0,0,0,0,0,0],#10          
[0,0,0,0,v,0,0,0,0,o,t,r,f,s,0,0,0,0,0,0,0,0,0,0,0,0,0,0,0,0,0,0,0,0,0,0,0,0,0,0,0,
0,0,0,0,0,0,0,0,0,0,0,0,0,0,0,0,0,0,0,0,0,0,0,0,0,0,0,0,0,0,0,0,0,0,0,0,0,0,0,0],#11
[0,0,0,0,f,0,0,0,0,0,o,t,r,v,s,0,0,0,0,0,0,0,0,0,0,0,0,0,0,0,0,0,0,0,0,0,0,0,0,0,0,
0,0,0,0,0,0,0,0,0,0,0,0,0,0,0,0,0,0,0,0,0,0,0,0,0,0,0,0,0,0,0,0,0,0,0,0,0,0,0,0],#12      
[0,0,0,0,r,0,0,0,0,0,0,o,t,f,v,s,0,0,0,0,0,0,0,0,0,0,0,0,0,0,0,0,0,0,0,0,0,0,0,0,0,
0,0,0,0,0,0,0,0,0,0,0,0,0,0,0,0,0,0,0,0,0,0,0,0,0,0,0,0,0,0,0,0,0,0,0,0,0,0,0,0],#13
[0,0,0,0,t,0,0,0,0,0,0,0,o,r,f,v,s,0,0,0,0,0,0,0,0,0,0,0,0,0,0,0,0,0,0,0,0,0,0,0,0,
0,0,0,0,0,0,0,0,0,0,0,0,0,0,0,0,0,0,0,0,0,0,0,0,0,0,0,0,0,0,0,0,0,0,0,0,0,0,0,0],#14
[0,0,0,0,o,0,0,0,0,0,0,0,0,t,r,f,v,0,0,0,0,0,0,0,0,0,0,0,0,0,0,0,0,0,0,s,0,0,0,0,0,
0,0,0,0,0,0,0,0,0,0,0,0,0,0,0,0,0,0,0,0,0,0,0,0,0,0,0,0,0,0,0,0,0,0,0,0,0,0,0,0],#15
[0,0,0,0,0,0,0,0,0,0,0,0,0,0,o,t,r,v,s,0,0,0,0,0,0,0,0,0,0,0,0,0,0,0,0,f,0,0,0,0,0,
0,0,0,0,0,0,0,0,0,0,0,0,0,0,0,0,0,0,0,0,0,0,0,0,0,0,0,0,0,0,0,0,0,0,0,0,0,0,0,0],#17
[0,0,0,0,0,0,0,0,0,0,0,0,0,0,0,o,t,f,v,s,0,0,0,0,0,0,0,0,0,0,0,0,0,0,0,r,0,0,0,0,0,
0,0,0,0,0,0,0,0,0,0,0,0,0,0,0,0,0,0,0,0,0,0,0,0,0,0,0,0,0,0,0,0,0,0,0,0,0,0,0,0],#18
[0,0,0,0,0,0,0,0,0,0,0,0,0,0,0,0,o,r,f,v,s,0,0,0,0,0,0,0,0,0,0,0,0,0,0,t,0,0,0,0,0,
0,0,0,0,0,0,0,0,0,0,0,0,0,0,0,0,0,0,0,0,0,0,0,0,0,0,0,0,0,0,0,0,0,0,0,0,0,0,0,0],#19
[0,0,0,0,0,0,0,0,0,0,0,0,0,0,0,0,0,t,r,f,v,s,0,0,0,0,0,0,0,0,0,0,0,0,0,o,0,0,0,0,0,
0,0,0,0,0,0,0,0,0,0,0,0,0,0,0,0,0,0,0,0,0,0,0,0,0,0,0,0,0,0,0,0,0,0,0,0,0,0,0,0],#20
[0,0,0,0,0,0,0,0,0,0,0,0,0,0,0,0,0,0,o,t,r,f,v,0,0,0,0,0,0,0,0,0,0,0,0,0,0,0,0,0,0,
0,0,0,0,0,0,0,0,0,0,0,0,0,0,0,0,0,0,0,0,0,0,0,0,0,0,0,0,s,0,0,0,0,0,0,0,0,0,0,0],#22
[0,0,0,0,0,0,0,0,0,0,0,0,0,0,0,0,0,0,0,o,t,r,f,s,0,0,0,0,0,0,0,0,0,0,0,0,0,0,0,0,0,
0,0,0,0,0,0,0,0,0,0,0,0,0,0,0,0,0,0,0,0,0,0,0,0,0,0,0,0,v,0,0,0,0,0,0,0,0,0,0,0],#23
[0,0,0,0,0,0,0,0,0,0,0,0,0,0,0,0,0,0,0,0,o,t,r,v,s,0,0,0,0,0,0,0,0,0,0,0,0,0,0,0,0,
0,0,0,0,0,0,0,0,0,0,0,0,0,0,0,0,0,0,0,0,0,0,0,0,0,0,0,0,f,0,0,0,0,0,0,0,0,0,0,0],#24
[0,0,0,0,0,0,0,0,0,0,0,0,0,0,0,0,0,0,0,0,0,o,t,f,v,s,0,0,0,0,0,0,0,0,0,0,0,0,0,0,0,
0,0,0,0,0,0,0,0,0,0,0,0,0,0,0,0,0,0,0,0,0,0,0,0,0,0,0,0,r,0,0,0,0,0,0,0,0,0,0,0],#25
[0,0,0,0,0,0,0,0,0,0,0,0,0,0,0,0,0,0,0,0,0,0,o,r,f,v,s,0,0,0,0,0,0,0,0,0,0,0,0,0,0,
0,0,0,0,0,0,0,0,0,0,0,0,0,0,0,0,0,0,0,0,0,0,0,0,0,0,0,0,t,0,0,0,0,0,0,0,0,0,0,0],#26
[0,0,0,0,0,0,0,0,0,0,0,0,0,0,0,0,0,0,0,0,0,0,0,t,r,f,v,s,0,0,0,0,0,0,0,0,0,0,0,0,0,
0,0,0,0,0,0,0,0,0,0,0,0,0,0,0,0,0,0,0,0,0,0,0,0,0,0,0,0,o,0,0,0,0,0,0,0,0,0,0,0],#27
[0,0,0,0,0,0,0,0,0,0,0,0,0,0,0,0,0,0,0,0,0,0,0,0,o,t,r,f,v,s,0,0,0,0,0,0,0,0,0,0,0,
0,0,0,0,0,0,0,0,0,0,0,0,0,0,0,0,0,0,0,0,0,0,0,0,0,0,0,0,0,0,0,0,0,0,0,0,0,0,0,0],#29           
[0,0,0,0,0,0,0,0,0,0,0,0,0,0,0,0,0,0,0,0,0,0,0,0,0,o,t,r,f,v,0,0,0,0,0,0,0,s,0,0,0,
0,0,0,0,0,0,0,0,0,0,0,0,0,0,0,0,0,0,0,0,0,0,0,0,0,0,0,0,0,0,0,0,0,0,0,0,0,0,0,0],#30
[0,0,0,0,0,0,0,0,0,0,0,0,0,0,0,0,0,0,0,0,0,0,0,0,0,0,o,t,r,f,s,0,0,0,0,0,0,v,0,0,0,
0,0,0,0,0,0,0,0,0,0,0,0,0,0,0,0,0,0,0,0,0,0,0,0,0,0,0,0,0,0,0,0,0,0,0,0,0,0,0,0],#31
[0,0,0,0,0,0,0,0,0,0,0,0,0,0,0,0,0,0,0,0,0,0,0,0,0,0,0,o,t,r,v,s,0,0,0,0,0,f,0,0,0,
0,0,0,0,0,0,0,0,0,0,0,0,0,0,0,0,0,0,0,0,0,0,0,0,0,0,0,0,0,0,0,0,0,0,0,0,0,0,0,0],#32
[0,0,0,0,0,0,0,0,0,0,0,0,0,0,0,0,0,0,0,0,0,0,0,0,0,0,0,0,o,t,f,v,s,0,0,0,0,r,0,0,0,
0,0,0,0,0,0,0,0,0,0,0,0,0,0,0,0,0,0,0,0,0,0,0,0,0,0,0,0,0,0,0,0,0,0,0,0,0,0,0,0],#33
[0,0,0,0,0,0,0,0,0,0,0,0,0,0,0,0,0,0,0,0,0,0,0,0,0,0,0,0,0,o,r,f,v,s,0,0,0,t,0,0,0,
0,0,0,0,0,0,0,0,0,0,0,0,0,0,0,0,0,0,0,0,0,0,0,0,0,0,0,0,0,0,0,0,0,0,0,0,0,0,0,0],#34
[0,0,0,0,0,0,0,0,0,0,0,0,0,0,0,0,0,0,0,0,0,0,0,0,0,0,0,0,0,0,t,r,f,v,s,0,0,o,0,0,0,
0,0,0,0,0,0,0,0,0,0,0,0,0,0,0,0,0,0,0,0,0,0,0,0,0,0,0,0,0,0,0,0,0,0,0,0,0,0,0,0],#35
[0,0,0,0,0,0,0,0,0,0,0,0,0,0,0,0,0,0,0,0,0,0,0,0,0,0,0,0,0,0,0,o,t,r,f,v,s,0,0,0,0,
0,0,0,0,0,0,0,0,0,0,0,0,0,0,0,0,0,0,0,0,0,0,0,0,0,0,0,0,0,0,0,0,0,0,0,0,0,0,0,0],#37
[0,0,0,0,0,0,0,0,0,0,0,0,0,0,0,0,0,0,0,0,0,0,0,0,0,0,0,0,0,0,0,0,o,t,r,f,v,s,0,0,0,
0,0,0,0,0,0,0,0,0,0,0,0,0,0,0,0,0,0,0,0,0,0,0,0,0,0,0,0,0,0,0,0,0,0,0,0,0,0,0,0],#38
[0,0,0,0,0,0,0,0,0,0,0,0,0,0,0,0,0,0,0,0,0,0,0,0,0,0,0,0,0,0,0,0,0,o,t,r,f,v,s,0,0,
0,0,0,0,0,0,0,0,0,0,0,0,0,0,0,0,0,0,0,0,0,0,0,0,0,0,0,0,0,0,0,0,0,0,0,0,0,0,0,0],#39
[0,0,0,0,0,0,0,0,0,0,0,0,0,0,0,0,0,0,0,0,0,0,0,0,0,0,0,0,0,0,0,0,0,0,o,t,r,f,v,s,0,
0,0,0,0,0,0,0,0,0,0,0,0,0,0,0,0,0,0,0,0,0,0,0,0,0,0,0,0,0,0,0,0,0,0,0,0,0,0,0,0],#40
[0,0,0,0,0,0,0,0,0,0,0,0,0,0,0,0,0,0,0,0,0,0,0,0,0,0,0,0,0,0,0,0,0,0,0,o,t,r,f,v,s,
0,0,0,0,0,0,0,0,0,0,0,0,0,0,0,0,0,0,0,0,0,0,0,0,0,0,0,0,0,0,0,0,0,0,0,0,0,0,0,0],#41
[0,0,0,0,0,0,0,0,0,0,0,0,0,0,0,0,0,0,0,0,0,s,0,0,0,0,0,0,0,0,0,0,0,0,0,0,o,t,r,f,v,
0,0,0,0,0,0,0,0,0,0,0,0,0,0,0,0,0,0,0,0,0,0,0,0,0,0,0,0,0,0,0,0,0,0,0,0,0,0,0,0],#42
[0,0,0,0,0,0,0,0,s,0,0,0,0,0,0,0,0,0,0,0,0,v,0,0,0,0,0,0,0,0,0,0,0,0,0,0,0,o,t,r,f,
0,0,0,0,0,0,0,0,0,0,0,0,0,0,0,0,0,0,0,0,0,0,0,0,0,0,0,0,0,0,0,0,0,0,0,0,0,0,0,0],#43
[0,0,0,0,0,0,0,0,v,0,0,0,0,0,0,0,0,0,0,0,0,f,0,0,0,0,0,0,0,0,0,0,0,0,0,0,0,0,o,t,r,
s,0,0,0,0,0,0,0,0,0,0,0,0,0,0,0,0,0,0,0,0,0,0,0,0,0,0,0,0,0,0,0,0,0,0,0,0,0,0,0],#44
[0,0,0,0,0,0,0,0,f,0,0,0,0,0,0,0,0,0,0,0,0,r,0,0,0,0,0,0,0,0,0,0,0,0,0,0,0,0,0,o,t,
v,0,0,0,0,0,0,0,0,0,0,0,0,s,0,0,0,0,0,0,0,0,0,0,0,0,0,0,0,0,0,0,0,0,0,0,0,0,0,0],#45
[0,0,0,0,0,0,0,0,r,0,0,0,0,0,0,0,0,0,0,0,0,t,0,0,0,0,0,0,0,0,0,0,0,0,0,0,0,0,0,0,o,
f,s,0,0,0,0,0,0,0,0,0,0,0,v,0,0,0,0,0,0,0,0,0,0,0,0,0,0,0,0,0,0,0,0,0,0,0,0,0,0],#46
[0,0,0,0,0,0,0,0,t,0,0,0,0,0,0,0,0,0,0,0,0,o,0,0,0,0,0,0,0,0,0,0,0,0,0,0,0,0,0,0,0,
r,v,s,0,0,0,0,0,0,0,0,0,0,f,0,0,0,0,0,0,0,0,0,0,0,0,0,0,0,0,0,0,0,0,0,0,0,0,0,0],#47
[0,0,0,0,0,0,0,0,0,0,0,0,0,0,0,0,0,0,0,0,0,0,0,0,0,0,0,0,0,0,0,0,0,0,0,0,0,0,0,0,0,
0,t,r+s,f,v,0,0,0,0,0,0,0,0,o,0,0,0,0,0,0,0,0,0,0,0,0,0,0,0,0,0,0,0,0,0,0,0,0,0,0],#50
[0,0,0,0,0,0,0,0,0,0,0,0,0,0,0,0,0,0,0,0,0,0,0,0,0,0,0,0,0,0,0,0,0,0,0,0,0,0,0,0,0,
0,0,o+f,t,r,v,s,0,0,0,0,0,0,0,0,0,0,0,0,0,0,0,0,0,0,0,0,0,0,0,0,0,0,0,0,0,0,0,0,0],#52
[0,0,0,0,0,0,0,0,0,0,0,0,0,0,0,0,0,0,0,0,0,0,0,0,0,0,0,0,0,0,0,0,0,0,0,0,0,0,0,0,0,
0,0,r,o,t,f,v,s,0,0,0,0,0,0,0,0,0,0,0,0,0,0,0,0,0,0,0,0,0,0,0,0,0,0,0,0,0,0,0,0],#53
[0,0,0,0,0,0,0,0,0,0,0,0,0,0,0,0,0,0,0,0,0,0,0,0,0,0,0,0,0,0,0,0,0,0,0,0,0,0,0,0,0,
0,0,t,0,o,r,f,v,s,0,0,0,0,0,0,0,0,0,0,0,0,0,0,0,0,0,0,0,0,0,0,0,0,0,0,0,0,0,0,0],#54
[0,0,0,0,0,0,0,0,0,0,0,0,0,0,0,0,0,0,0,0,0,0,0,0,0,0,0,0,0,0,0,0,0,0,0,0,0,0,0,0,0,
0,0,o,0,0,t,r,f,v,s,0,0,0,0,0,0,0,0,0,0,0,0,0,0,0,0,0,0,0,0,0,0,0,0,0,0,0,0,0,0],#55
[0,0,0,0,0,0,0,0,0,0,0,0,0,0,0,v,0,0,0,0,0,0,0,0,0,0,0,0,0,0,0,0,0,0,0,0,0,0,0,0,0,
0,0,0,0,0,0,o,t,r,f,s,0,0,0,0,0,0,0,0,0,0,0,0,0,0,0,0,0,0,0,0,0,0,0,0,0,0,0,0,0],#57
[0,0,0,0,0,0,0,0,0,0,0,0,0,0,0,f,0,0,0,0,0,0,0,0,0,0,0,0,0,0,0,0,0,0,0,0,0,0,0,0,0,
0,0,0,0,0,0,0,o,t+s,r,v,0,0,0,0,0,0,0,0,0,0,0,0,0,0,0,0,0,0,0,0,0,0,0,0,0,0,0,0,0],#58
[0,0,0,0,0,0,0,0,0,0,0,0,0,0,0,r,0,0,0,0,0,0,0,0,0,0,0,0,0,0,0,0,0,0,0,0,0,0,0,0,0,
0,0,0,0,0,0,0,0,o+v,t,f,s,0,0,0,0,0,0,0,0,0,0,0,0,0,0,0,0,0,0,0,0,0,0,0,0,0,0,0,0],#59
[0,0,0,0,0,0,0,0,0,0,0,0,0,0,0,t,0,0,0,0,0,0,0,0,0,0,0,0,0,0,0,0,0,0,0,0,0,0,0,0,0,
0,0,0,0,0,0,0,0,f,o,r,v,s,0,0,0,0,0,0,0,0,0,0,0,0,0,0,0,0,0,0,0,0,0,0,0,0,0,0,0],#60
[0,0,0,0,0,0,0,0,0,0,0,0,0,0,0,o,0,0,0,0,0,0,0,0,0,0,0,0,0,0,0,0,0,0,0,0,0,0,0,0,0,
0,0,0,0,0,0,0,0,r,0,t,f,v,s,0,0,0,0,0,0,0,0,0,0,0,0,0,0,0,0,0,0,0,0,0,0,0,0,0,0],#61
[0,0,0,0,0,0,0,0,0,0,0,0,0,0,0,0,0,0,0,0,0,0,0,0,0,0,0,0,0,0,0,0,0,0,0,0,0,0,0,0,0,
0,0,0,0,0,0,0,0,o,0,0,t,r,f,v,s,0,0,0,0,0,0,0,0,0,0,0,0,0,0,0,0,0,0,0,0,0,0,0,0],#63
[0,0,0,0,0,0,0,0,0,0,0,0,0,0,0,0,0,0,0,0,0,0,0,0,0,0,0,0,0,0,0,0,0,0,0,0,0,0,0,0,0,
0,0,0,0,0,0,0,0,0,0,0,0,o,t,r,f,v,0,0,0,0,0,0,0,0,0,0,0,0,0,0,0,0,0,s,0,0,0,0,0],#65
[0,0,0,0,0,0,0,0,0,0,0,0,0,0,0,0,0,0,0,0,0,0,0,0,0,0,0,0,0,0,0,0,0,0,0,0,0,0,0,0,0,
0,0,0,0,0,0,0,0,0,0,0,0,0,o,t,r,f,s,0,0,0,0,0,0,0,0,0,0,0,0,0,0,0,0,v,0,0,0,0,0],#66
[0,0,0,0,0,0,0,0,0,0,0,0,0,0,0,0,0,0,0,0,0,0,0,0,0,0,0,0,0,0,0,0,0,0,0,0,0,0,0,0,0,
0,0,0,0,0,0,0,0,0,0,0,0,0,0,o,t,r,v,s,0,0,0,0,0,0,0,0,0,0,0,0,0,0,0,f,0,0,0,0,0],#67
[0,0,0,0,0,0,0,0,0,0,0,0,0,0,0,0,0,0,0,0,0,0,0,0,0,0,0,0,0,0,0,0,0,0,0,0,0,0,0,0,0,
0,0,0,0,0,0,0,0,0,0,0,0,0,0,0,o,t,f,v,s,0,0,0,0,0,0,0,0,0,0,0,0,0,0,r,0,0,0,0,0],#68
[0,0,0,0,0,0,0,0,0,0,0,0,0,0,0,0,0,0,0,0,0,0,0,0,0,0,0,0,0,0,0,0,0,0,0,0,0,0,0,0,0,
0,0,0,0,0,0,0,0,0,0,0,0,0,0,0,0,o,r,f,v,s,0,0,0,0,0,0,0,0,0,0,0,0,0,t,0,0,0,0,0],#69
[0,0,0,0,0,0,0,0,0,0,0,0,0,0,0,0,0,0,0,0,0,0,0,0,0,0,0,0,0,0,0,0,0,0,0,0,0,0,0,0,0,
0,0,0,0,0,0,0,0,0,0,0,0,0,0,0,0,0,t,r,f,v,s,0,0,0,0,0,0,0,0,0,0,0,0,o,0,0,0,0,0],#70
[0,0,0,0,0,0,0,0,0,0,0,0,0,0,0,0,0,0,0,0,0,0,0,0,0,0,0,0,0,0,0,0,0,0,0,0,0,0,0,0,0,
0,0,0,0,0,0,0,0,0,0,0,0,0,0,0,0,0,0,o,t,r,f,v,s,0,0,0,0,0,0,0,0,0,0,0,0,0,0,0,0],#72
[0,0,0,0,0,0,0,0,0,0,0,0,0,0,0,0,0,0,0,0,0,0,0,0,0,0,0,0,0,0,0,0,0,0,0,0,0,0,0,0,0,
0,0,0,0,0,0,0,0,0,0,0,0,0,0,0,0,0,0,0,o,t,r,f,v,s,0,0,0,0,0,0,0,0,0,0,0,0,0,0,0],#73
[0,0,0,0,0,0,0,0,0,0,0,0,0,0,0,0,0,0,0,0,0,0,0,0,0,0,0,0,0,0,0,0,0,0,0,0,0,0,0,0,0,
0,0,0,0,0,0,0,0,0,0,0,0,0,0,0,0,0,0,0,0,o,t,r,f,v,0,0,0,0,0,0,0,0,0,0,0,0,0,0,0],#74
[0,0,0,0,0,0,0,0,0,0,0,0,0,0,0,0,0,0,0,0,0,0,0,0,0,0,0,0,0,0,0,0,0,0,0,0,0,0,0,0,0,
0,0,0,0,0,0,0,0,0,0,0,0,0,0,0,0,0,0,0,0,0,o,t,r,f,s,0,0,0,0,0,0,0,0,0,0,0,0,0,0],#75
[0,0,0,0,0,0,0,0,0,0,0,0,0,0,0,0,0,0,0,0,0,0,0,0,0,0,0,0,0,0,0,0,0,0,0,0,0,0,0,0,0,
0,0,0,0,0,0,0,0,0,0,0,0,0,0,0,0,0,0,0,0,0,0,o,t,r,v,s,0,0,0,0,0,0,0,0,0,0,0,0,0],#76
[0,0,0,0,0,0,0,0,0,0,0,0,0,0,0,0,0,0,0,0,0,0,0,0,0,0,0,0,0,0,0,0,0,0,0,0,0,0,0,0,0,
0,0,0,0,0,0,0,0,0,0,0,0,0,0,0,0,0,0,0,0,0,0,0,o,t,f,v,s,0,0,0,0,0,0,0,0,0,0,0,0],#77
[0,0,0,0,0,0,0,0,0,0,0,0,0,0,0,0,0,0,0,0,0,0,0,0,0,0,0,0,0,0,0,0,0,0,0,0,0,0,0,0,0,
0,0,0,0,0,0,0,0,0,0,0,0,0,0,0,0,0,0,0,0,0,0,0,0,o,r,f,v,s,0,0,0,0,0,0,0,0,0,0,0],#78
[0,0,0,0,0,0,0,0,0,0,0,0,0,0,0,0,0,0,0,0,0,0,0,0,0,0,0,0,0,0,0,0,0,0,0,0,0,0,0,0,0,
0,0,0,0,0,0,0,0,0,0,0,0,0,0,0,0,0,0,0,0,0,0,0,0,0,t,r,f,v,s,0,0,0,0,0,0,0,0,0,0],#79
[0,0,0,0,0,0,0,0,0,0,0,0,0,0,0,0,0,0,0,s,0,0,0,0,0,0,0,0,0,0,0,0,0,0,0,0,0,0,0,0,0,
0,0,0,0,0,0,0,0,0,0,0,0,0,0,0,0,0,0,0,0,0,0,0,0,0,0,o,t,r,f,v,0,0,0,0,0,0,0,0,0],#81
[0,0,0,0,0,0,0,0,0,0,0,0,0,0,0,0,0,0,0,v,0,0,0,0,0,0,0,0,0,0,0,0,0,0,0,0,0,0,0,0,0,
0,0,0,0,0,0,0,0,0,0,0,0,0,0,0,0,0,0,0,0,0,0,0,0,0,0,0,o,t,r,f,s,0,0,0,0,0,0,0,0],#82
[0,0,0,0,0,0,0,0,0,0,0,0,0,0,0,0,0,0,0,f,0,0,0,0,0,0,0,0,0,0,0,0,0,0,0,0,0,0,0,0,0,
0,0,0,0,0,0,0,0,0,0,0,0,0,0,0,0,0,0,0,0,0,0,0,0,0,0,0,0,o,t,r,v,s,0,0,0,0,0,0,0],#83
[0,0,0,0,0,0,0,0,0,0,0,0,0,0,0,0,0,0,0,r,0,0,0,0,0,0,0,0,0,0,0,0,0,0,0,0,0,0,0,0,0,
0,0,0,0,0,0,0,0,0,0,0,0,0,0,0,0,0,0,0,0,0,0,0,0,0,0,0,0,0,o,t,f,v,s,0,0,0,0,0,0],#84
[0,0,0,0,0,0,0,0,0,0,0,0,0,0,0,0,0,0,0,t,0,0,0,0,0,0,0,0,0,0,0,0,0,0,0,0,0,0,0,0,0,
0,0,0,0,0,0,0,0,0,0,0,0,0,0,0,0,0,0,0,0,0,0,0,0,0,0,0,0,0,0,o,r,f,v,s,0,0,0,0,0],#85
[0,0,0,0,0,0,0,0,0,0,0,0,0,0,0,0,0,0,0,o,0,0,0,0,0,0,0,0,0,0,0,0,0,0,0,0,0,0,0,0,0,
0,0,0,0,0,0,0,0,0,0,0,0,0,0,0,0,0,0,0,0,0,0,0,0,0,0,0,0,0,0,0,t,r,f,v,s,0,0,0,0],#86
[0,0,0,0,0,0,0,0,0,0,0,0,0,0,0,0,0,0,0,0,0,0,0,0,0,0,0,0,0,0,0,0,0,0,0,0,0,0,0,0,0,
0,0,0,0,0,0,0,0,0,0,0,0,0,0,0,0,0,0,v,0,0,0,0,0,0,0,0,0,0,0,0,0,o,t,r,f,s,0,0,0],#88
[0,0,0,0,0,0,0,0,0,0,0,0,0,0,0,0,0,0,0,0,0,0,0,0,0,0,0,0,0,0,0,0,0,0,0,0,0,0,0,0,0,
0,0,0,0,0,0,0,0,0,0,0,0,0,0,0,0,0,0,f,0,s,0,0,0,0,0,0,0,0,0,0,0,0,o,t,r,v,0,0,0],#89
[0,0,0,0,0,0,0,0,0,0,0,0,0,0,0,0,0,0,0,0,0,0,0,0,0,0,0,0,0,0,0,0,0,0,0,0,0,0,0,0,0,
0,0,0,0,0,0,0,0,0,0,0,0,0,0,0,0,0,0,r,0,v,0,0,0,0,0,0,0,0,0,0,0,0,0,o,t,f,s,0,0],#90
[0,0,0,0,0,0,0,0,0,0,0,0,0,0,0,0,0,0,0,0,0,0,0,0,0,0,0,0,0,0,0,0,0,0,0,0,0,0,0,0,0,
0,0,0,0,0,0,0,0,0,0,0,0,0,0,0,0,0,0,t,0,f,0,0,0,0,0,0,0,0,0,0,0,0,0,0,o,r,v,s,0],#91
[0,0,0,0,0,0,0,0,0,0,0,0,0,0,0,0,0,0,0,0,0,0,0,0,0,0,0,0,0,0,0,0,0,0,0,0,0,0,0,0,0,
0,0,0,0,0,0,0,0,0,0,0,0,0,0,0,0,0,0,o,0,r,0,0,s,0,0,0,0,0,0,0,0,0,0,0,0,t,f,v,0],#92
[0,0,0,0,0,0,0,0,0,0,0,0,0,0,0,0,0,0,0,0,0,0,0,0,0,0,0,0,0,0,0,0,0,0,0,0,0,0,0,0,0,
0,0,0,0,0,0,0,0,0,0,0,0,0,0,0,0,0,0,0,0,o,0,0,f,0,0,0,0,0,0,0,0,0,0,0,0,0,t,r,v],#94
[0,0,0,0,0,0,0,0,0,0,0,0,0,0,0,0,0,0,0,0,0,0,0,0,0,0,0,0,0,0,0,0,0,0,0,0,0,0,0,0,0,
0,0,0,0,0,0,0,0,0,0,0,0,0,0,0,0,0,0,0,0,0,0,0,t,0,0,0,0,0,0,0,0,0,0,0,0,0,v+s,o,r],#96
[0,0,0,0,0,0,0,0,0,0,0,0,0,0,0,0,0,0,0,0,0,0,0,0,0,0,0,0,0,0,0,0,0,0,0,0,0,0,0,0,0,
0,0,0,0,0,0,0,0,0,0,0,0,0,0,0,0,0,0,0,0,0,0,0,o,0,0,0,0,0,0,0,0,0,0,0,0,0,0,f+v+s,t],#97
[0,0,0,0,0,0,0,0,0,0,0,0,0,0,0,0,0,0,0,0,0,0,0,0,0,0,0,0,0,0,0,0,0,0,0,0,0,0,0,0,0,
0,0,0,0,0,0,0,0,0,0,0,0,0,0,0,0,0,0,0,0,0,0,0,0,0,0,0,0,0,0,0,0,0,0,0,0,0,0,0,t+r+f+v+s]])#99

if o + t + r + f + v + s != 1:
    print("error")
(spos*((matrix.identity(81) - T).inverse())*one).n()

\end{lstlisting}
\end{center}

\newpage

\lstset{
    language=python,
    tabsize=2,
    %frame=lines,
    frame=shadowbox,
    rulesepcolor=\color{gray},
    xleftmargin=-50pt,
    xrightmargin=-50pt,
    framexleftmargin=15pt,
    keywordstyle=\color{blue},
    commentstyle=\color{60color},
    stringstyle=\color{red},
    numbers=left,
    numberstyle=\tiny,
    numbersep=5pt,
    breaklines=true,
    %postbreak=\mbox{\textcolor{red}{$\hookrightarrow$}\space},
    showstringspaces=false,
    basicstyle=\footnotesize,
    emph={str},emphstyle={\color{magenta}},
    caption={Python Code: Simulate All Start Positions Over All $P(\delta)=1$}, captionpos=t
    }
    
\label{py2}
\begin{lstlisting}
import random


def Welcome():
    line = "Weighted Number, Start Pos, Result, Moves, Pos"
    with open("sposData.txt",'w') as file:
        file.write(line)
    print("Welcome to the Chutes & Ladders simulation!")
    Start()

def Start():
    #initialize vars
    nums = [1,2,3,4,5,6]
    weight_num = 0
    w=[0]

    for r in range (1,7,1):
        w=[0]
        nums = [1,2,3,4,5,6]
        weight_num = r
        nums.remove(int(weight_num))
        for i in range(1,7,1):
            if nums.count(i) > 0:
                if i > weight_num:
                    w.append(10000000)
                else:
                    w.append(0)
            else:
                w.append(10000000)
        print("For", r,":", w)
        for spos in range (0,100,1):
            Game(w, spos, weight_num)

    print("Done!")

def dist(w):
    broll = random.randint(1,10000000)
    if broll <= w[1]:
        return int(1)
    elif broll <= w[2]:
        return int(2)
    elif broll <= w[3]:
        return int(3)
    elif broll <= w[4]:
        return int(4)
    elif broll <= w[5]:
        return int(5)
    elif broll <= w[6]:
        return int(6)





def Game(w, spos, weight_num):
    game_num = 1
    long_counter = 0
    last_5 = []
    loop = ""
    # Set up the board #
    #Pos:       1   2  3  4   5  6  7  8  9  10 11 12 13 14 15 16 17 18 19 20  21 22 23 24 25 26 27  28 29 30 31 32 33 34 35  36 37 38 39 40 41 42 43 44 45 46 47  48  49 50
    board = [0, 38, 0, 0, 14, 0, 0, 0, 0, 31, 0, 0, 0, 0, 0, 0, 6, 0, 0, 0, 0, 42, 0, 0, 0, 0, 0, 0, 84, 0, 0, 0, 0, 0, 0, 0, 44, 0, 0, 0, 0, 0, 0, 0, 0, 0, 0, 0, 26, 11, 0,
                67, 0, 0, 0, 0, 53, 0, 0, 0, 0, 0, 19, 0, 60, 0, 0, 0, 0, 0, 0, 91, 0, 0, 0, 0, 0, 0, 0, 0, 100, 0, 0, 0, 0, 0, 0, 24, 0, 0, 0, 0, 0, 73, 0, 75, 0, 0, 78, 0, 0]
    #Pos:       51 52 53 54 55  56 57 58 59 60 61  62 63  64 65 66 67 68 69 70  71 72 73 74 75 76 77 78 79   80 81 82 83 84 85 86  87 88 89 90 91 92  93 94  95 96 97  98 99 100

    #Initialize move counter
    move_number = 0


    # Initialize player position
    player_pos = spos
    if board[spos] > 0:
        player_pos = board[spos]
        if player_pos == 100:
            # WIINER WINNER CHICEKN DINNER
            gameOn = False
            file_write(str(str(weight_num) + "," + str(spos) + "," + "Win," + str(move_number) + "," + str(player_pos)))
        else:
            gameOn = True
    else:
        gameOn = True







    while gameOn == True:
        move_number += 1

        #Roll
        roll = dist(w)





        # Move the player
        player_pos += roll


        #Check new location
        #Finale Loop
        if player_pos > 100:
            player_pos -= roll
            gameOn = False
            file_write(str(str(weight_num) + "," + str(spos) + "," + "Finale Loop," + "Inf," + str(player_pos)))
        else:
            if board[player_pos] != 0:
                player_pos = board[player_pos]


        if player_pos == 100:
            # WIINER WINNER CHICEKN DINNER
            gameOn = False
            file_write(str(str(weight_num) + "," + str(spos) + "," + "Win," + str(move_number) + "," + str(player_pos)))
        if len(last_5) < 6:
            last_5.append(player_pos)
        else:
            last_5.pop(0)
            last_5.append(player_pos)

        if move_number > 500:
            loop = ' '.join([str(elem) for elem in last_5])
            gameOn = False
            file_write(str(str(weight_num) + "," + str(spos) + "," + "C&L Loop," + "Inf," + loop))






def file_write(line):
    with open("sposData.txt", 'a') as file:
        file.write('\n')
        file.write(str(line))






Welcome()
\end{lstlisting}

\newpage

\lstset{
    language=python,
    tabsize=2,
    %frame=lines,
    frame=shadowbox,
    rulesepcolor=\color{gray},
    xleftmargin=-50pt,
    xrightmargin=-50pt,
    framexleftmargin=15pt,
    keywordstyle=\color{blue},
    commentstyle=\color{60color},
    stringstyle=\color{red},
    numbers=left,
    numberstyle=\tiny,
    numbersep=5pt,
    breaklines=true,
    %postbreak=\mbox{\textcolor{red}{$\hookrightarrow$}\space},
    showstringspaces=false,
    basicstyle=\footnotesize,
    emph={str},emphstyle={\color{magenta}},
    caption={Python Code: Simulate Game With Coin Strategies}, captionpos=t
    }
    
\label{py3}
\begin{lstlisting}
import random
import math


def Welcome():
    print("Welcome to the Chutes & Ladders simulation!")
    Start()

def Start():
    #initialize vars
    tbt_bool = True
    sum_bool = False
    weight_bool = False
    total_games = ""
    strat = int(0)

    #input
    tbt_YN = input("Do you want a turn by turn readout? (Only able to simulate a game at a time) ")
    if tbt_YN == "Y" or tbt_YN == "y":
        tbt_bool = True
    elif tbt_YN == "N" or tbt_YN == "n":
        tbt_bool = False
    else:
        print("Only 'y' or 'n'.")
        Start()
    if tbt_bool == False:
        sum_YN = input("Do you want a summary readout for each game? ")
        if sum_YN == "Y" or sum_YN == "y":
            sum_bool = True
        elif sum_YN == "N" or sum_YN == "n":
            sum_bool = False
        else:
            print("Only 'y' or 'n'.")
            Start()

    if tbt_bool == True:
        total_games = 1
    else:
        total_games = input("Number of games to simulate: ")
        try:
            int(total_games)
        except:
            print("You must enter a positive integer.")
            Start()
        if int(total_games) <= 0:
            print("You must enter a positive integer.")
            Start()




    Game(tbt_bool,sum_bool,int(total_games),strat)



def my_round(x):
    frac = x - math.floor(x)
    if frac < 0.5:
        return math.floor(x)
    else:
        return math.ceil(x)


def strategy(strat, board, player_pos):
    flips = [-1,1]
    #If player landed on square-100, do not flip. Independent of stratm, except strat -1 (random).
    if strat != -1:
        if player_pos == 100:
            return 0

#Strat -1 is to randomly flip
    if strat == -1:
        if my_round(random.randint(0,100)/100) == 0:
            return flips[my_round(random.randint(0,100)/100)]
        else:
            return 0
#Strat 0 is the base game, never flip coin
    elif strat == 0:
        return 0
#Strat 1 os to flip randomly, except never flip on square-100
    elif strat == 1:
        if my_round(random.randint(0,100)/100) == 0:
            return flips[my_round(random.randint(0,100)/100)]
        else:
            return 0
#Strat 2 is to always flip coin
    elif strat == 2:
        return flips[my_round(random.randint(0,100)/100)] #Always flip coin
#Strat 3 is to always flip unless at base of ladder
    elif strat == 3:
        if player_pos > 100:
            return flips[my_round(random.randint(0,100)/100)] #Flip coin if they overshot square-100
        #If on top of slide:
        elif board[player_pos] < player_pos: #If top of chute or not base of ladder
            return flips[my_round(random.randint(0,100)/100)] #Flip coin
        #If at base of ladder:
        else: #Else = at base of ladder
            return 0 #Do not flip coin
#Strat 4 is to only flip if land at top of slide
    elif strat == 4:
        if player_pos > 100:
            return flips[my_round(random.randint(0,100)/100)] #Flip coin if they overshot square-100
        #If on top of slide:
        elif board[player_pos] < player_pos and board[player_pos] != 0:
            return flips[my_round(random.randint(0,100)/100)] #Flip coin
        #If at base of ladder:
        elif board[player_pos] > player_pos:
            return 0 #Do not flip coin
        else:
            return 0 #Do not flip coin
#Strat 5 is to only flip if at top of slide or 1 square from a ladder.
    elif strat == 5:
        if player_pos > 100:
            return flips[my_round(random.randint(0,100)/100)] #Flip coin if they overshot square-100
        #If on top of slide:
        elif board[player_pos] < player_pos and board[player_pos] != 0:
            return flips[my_round(random.randint(0,100)/100)] #Flip coin
        #If at base of ladder:
        elif board[player_pos] > player_pos:
            return 0 #Do not flip coin
        elif board[player_pos - 1] > player_pos or board[player_pos + 1] > player_pos: #If within on square of ladder
            return flips[my_round(random.randint(0,100)/100)]
        else:
            return 0
#Strat 6 is to never flip when at base of ladder or one square from top of slide. Otherwise, always flip.
    elif strat == 6:
        if player_pos > 100:
            return flips[my_round(random.randint(0,100)/100)] #Flip coin if they overshot square-100
        #If on top of slide:
        elif board[player_pos] < player_pos and board[player_pos] != 0:
            return flips[my_round(random.randint(0,100)/100)] #Flip coin
        #If at base of ladder:
        elif board[player_pos] > player_pos:
            return 0 #Do not flip coin
        elif board[player_pos - 1] < player_pos and board[player_pos - 1] != 0: #If tails would lead to slide
            return 0 #Do not flip coin
        elif board[player_pos + 1] < player_pos and board[player_pos + 1] != 0: #If heads would lead to slide
            return 0
        else:
            return flips[my_round(random.randint(0,100)/100)] #Flip coin





def Game(display, summary, total_games,strat):
    game_num = 1
    notified = False
    long_counter = 0
    # Set up the boar#
    #Pos:       1   2  3  4   5  6  7  8  9  10 11 12 13 14 15 16 17 18 19 20  21 22 23 24 25 26 27  28 29 30 31 32 33 34 35  36 37 38 39 40 41 42 43 44 45 46 47  48  49 50
    board = [0, 38, 0, 0, 14, 0, 0, 0, 0, 31, 0, 0, 0, 0, 0, 0, 6, 0, 0, 0, 0, 42, 0, 0, 0, 0, 0, 0, 84, 0, 0, 0, 0, 0, 0, 0, 44, 0, 0, 0, 0, 0, 0, 0, 0, 0, 0, 0, 26, 11, 0,
             67, 0, 0, 0, 0, 53, 0, 0, 0, 0, 0, 19, 0, 60, 0, 0, 0, 0, 0, 0, 91, 0, 0, 0, 0, 0, 0, 0, 0, 100, 0, 0, 0, 0, 0, 0, 24, 0, 0, 0, 0, 0, 73, 0, 75, 0, 0, 78, 0, 0]
    #Pos:       51 52 53 54 55  56 57 58 59 60 61  62 63  64 65 66 67 68 69 70  71 72 73 74 75 76 77 78 79   80 81 82 83 84 85 86  87 88 89 90 91 92  93 94  95 96 97  98 99 100


    #Base board with only ladders:
    #Pos:       1   2  3  4   5  6  7  8  9  10 11 12 13 14 15 16 17 18 19 20  21 22 23 24 25 26 27  28 29 30 31 32 33 34 35  36 37 38 39 40 41 42 43 44 45 46 47 48 49 50
    #board = [0, 38, 0, 0, 14, 0, 0, 0, 0, 31, 0, 0, 0, 0, 0, 0, 0, 0, 0, 0, 0, 42, 0, 0, 0, 0, 0, 0, 84, 0, 0, 0, 0, 0, 0, 0, 44, 0, 0, 0, 0, 0, 0, 0, 0, 0, 0, 0, 0, 0, 0,
                #67, 0, 0, 0, 0, 0, 0, 0, 0, 0, 0, 0, 0, 0, 0, 0, 0, 0, 0, 0, 91, 0, 0, 0, 0, 0, 0, 0, 0, 100, 0, 0, 0, 0, 0, 0, 0, 0, 0, 0, 0, 0, 0, 0, 0, 0, 0, 0, 0, 0]
    #Pos:       51 52 53 54 55 56 57 58 59 60 61 62 63 64 65 66 67 68 69 70  71 72 73 74 75 76 77 78 79   80 81 82 83 84 85 86 87 88 89 90 91 92 93 94 95 96 97 98 99 100

    # Initialize player position
    player_pos = 0

    #Initialize move counter
    move_number = 0
    #Initialize game status
    gameOn = True

    #Initialize total move record
    moves_per_game = []
    #Initialize total moves
    total_moves = 0
    #Initialize coin flip result
    flip_res = 0




    while gameOn == True:
        move_number += 1
        #1st Game
        if game_num == 1 and move_number == 1 and summary == True:
            print("Game #1!")
        # Roll the dice
        if display == True:
            print("Move #",move_number,":", sep='')




        #Roll
        roll = random.randint(1,6)







        # Move the player
        player_pos += roll


        #Check new location
        if player_pos > 101:
            player_pos -= roll
            if display == True:
                print("You have far overshot square-100. Roll again.")
                print("You have returned to square", player_pos)
        else:
            flip_res = strategy(strat, board, player_pos)
            if display == True:
                print("You rolled a ", roll, " to square-",player_pos, end='. ', sep='')
                if player_pos <= 100:
                    if board[player_pos] != 0:
                        print("This square is a", "chute" if board[player_pos] < player_pos else "ladder!")
                if flip_res == 0:
                    print("You did not flip the coin")
                else:
                    print("You flipped the coin and got", "heads!" if flip_res == 1 else "tails!")
            player_pos += flip_res


        #Check location after flip
        if player_pos > 100:
            player_pos = player_pos - (flip_res + roll)
            if display == True:
                print("You overshot square-100. Roll again.")
        if board[player_pos] != 0:
            if display == True:
                print("You landed on a", "chute" if board[player_pos] < player_pos else "ladder!")
            player_pos = board[player_pos]
        if display == True:
            print("You are now on square", player_pos)






        if player_pos == 100:
            # WIINER WINNER CHICEKN DINNER
            moves_per_game.append(move_number)
            if display == True:
                print("You won!")
            if summary == True:
                print("Total number of moves:", move_number)
            if game_num == total_games:
                gameOn = False
                for i in range(total_games):
                    total_moves += moves_per_game[i-1]
                print("___________________________________________")
                print("Simulation comeplete. Played ", game_num, " game(s). Total moves: ", total_moves, ". Average moves per game: ", total_moves/total_games,".", sep='')
                print("Number of games over 1k moves", long_counter)
            else:
                if notified == True:
                    #print("Game", game_num, "took", move_number, "moves.")
                    long_counter += 1
                game_num += 1
                move_number = 0
                player_pos = 0
                notified = False
                if summary == True:
                    print("___________________________________________")
                    print("NEW GAME! This is now game #",game_num,".",sep='')









Welcome()
\end{lstlisting}

\end{document}